\documentclass[11pt]{article}
\usepackage{amsmath}
\usepackage{amsfonts}
\usepackage{amssymb}
\usepackage{graphicx}
\usepackage{color}

\numberwithin{equation}{section}
\definecolor{mycolorred}{rgb}{1, 0, 0}

\newtheorem{theorem}{Theorem}[section]

\newtheorem{assumption}[theorem]{Assumption}

\newtheorem{lemma}[theorem]{Lemma}

\newtheorem{proposition}[theorem]{Proposition}
\newtheorem{remark}[theorem]{Remark}

\def\<{\langle}
\def\>{\rangle}

\def\R{{\mathbb R}}

\begin{document}

\title{Transfer of regularity for Markov semigroups}
\author{\textsc{Vlad Bally}\thanks{
LAMA (UMR CNRS, UPEMLV, UPEC), MathRisk INRIA, Universit\'e Paris-Est - 
\texttt{vlad.bally@u-pem.fr}} \\
\textsc{Lucia Caramellino}\thanks{
Dipartimento di Matematica, Universit\`a di Roma Tor Vergata, and
INdAM-GNAMPA - \texttt{caramell@mat.uniroma2.it}}}
\date{}
\maketitle

\begin{abstract}
\noindent{\parindent0pt We study the regularity of a Markov semigroup $
(P_t)_{t>0}$, that is, when $P_t(x,dy)=p_t(x,y)dy$ for a suitable smooth 
function $p_t(x,y)$. This is done by transferring the regularity from an 
approximating Markov semigroup sequence $(P^n_t)_{t>0}$, $n\in{\mathbb{N}}$,
whose associated densities $p^n_t(x,y)$ are smooth and can blow up as $
n\to\infty$. We use an interpolation type result and we show that if there 
exists a good equilibrium between the blow-up and the speed of convergence, 
then $P_{t}(x,dy)=p_{t}(x,y)dy$ and $p_{t}$ has some regularity properties.  
}
\end{abstract}

\tableofcontents


\medskip


\medskip\noindent
{\bf Acknowledgments.}
This research is partly funded by the B\'ezout Labex, funded by ANR, reference ANR-10-LABX-58.
L.C. also acknowledges the MIUR Excellence
Department Project  awarded to the Department of Mathematics, University of
Rome Tor Vergata. 
The Beyond Borders Project ``Asymptotic Methods in Probability'' is acknowledged by both authors.

\section{Introduction}

In this paper we study Markov semigroups, that is, positive semigroups $%
(P_t)_{t\geq 0}$, 
such that $P_t1=1$. The link with Markov processes is given by a family $%
P_t(x,dy)$, $t\geq 0$, $x\in{\mathbb{R}}^d$, of transition probability
measures in ${\mathbb{R}}^d$ such that 
\begin{equation*}
P_{t}f(x)=\int_{{\mathbb{R}}^d}f(y)P_t(x,dy),\quad t\geq 0.
\end{equation*}
We study here the regularity of $P_t$, which is the property $%
P_t(x,dy)=p_t(x,y)dy$, $t>0$, for a suitable smooth function $p_t(x,y)$, by
transferring the regularity from an approximating Markov semigroup sequence $%
(P^n_t)_{t\geq 0}$, $n\in{\mathbb{N}}$.

Hereafter we assume that the domain of the Markov semigroup $(P_t)_{t\geq 0}$
contains the Schwartz space ${\mathcal{S}({\mathbb{R}}^d)}$ of the $C^\infty(%
{\mathbb{R}}^d)$ functions all of whose derivatives are rapidly decreasing.
We assume that the semigroup is strongly continuous in its domain and
we call $L$ its infinitesimal generator. We suppose also that the domain of $L
$ contains ${\mathcal{S}({\mathbb{R}}^d)}$ and for every $f\in {\ \mathcal{S}%
({\mathbb{R}}^d)}$, $P_tf\in {\mathcal{S}({\mathbb{R}}^d)}$, $t\geq 0$, and $%
Lf\in {\mathcal{S}({\mathbb{R}}^d)}$.

Let $(P^n_t)_{t\geq 0}$
, $n\in {\ \mathbb{N}}$, be a sequence of Markov semigroups: 
\begin{equation*}
P^n_{t}f(x)=\int_{{\mathbb{R}}^d}f(y)P^n_t(x,dy),\quad t\geq 0, n\in{\mathbb{%
\ N}}.
\end{equation*}
For every $n$, we assume that $(P^n_t)_{t\geq 0}$ satisfies the same
properties as $(P_t)_{t\geq 0}$: ${\mathcal{S}({\ \mathbb{R}}^{d})}$ is
included in the domain of $P_{t}^{n}$ and if $f\in {\mathcal{S}({\mathbb{R}}%
^d)}$ then $P_tf\in\mathcal{S}({\mathbb{R}}^d)$, $t\geq 0$; $(P^n_t)_{t\geq
0}$ is strongly continuous in its domain; the domain of its infinitesimal
operator $L_{n}$ contains ${\mathcal{S}({\mathbb{R}}^{d})}$ and $L_nf\in {%
\mathcal{S}({\mathbb{R}}^d)}$ if $f\in {\mathcal{S}({\mathbb{R}} ^d)}$.

Classical results (Trotter Kato theorem, see e.g. \cite{EK}) assert that, as 
$n\to\infty$, if $L_{n}\rightarrow L$ then $P_{t}^{n}\rightarrow P_{t}.$ The
problem that we address in this paper is the following. We suppose that $%
P_{t}^{n}$ has the regularity (density) property $%
P_{t}^{n}(x,dy)=p_{t}^{n}(x,y)dy$ with $p_{t}^{n}\in C^{\infty }({\mathbb{R}}%
^{d}\times {\mathbb{R}}^{d})$ and we ask under which hypotheses this
property is inherited by the limit semigroup $(P_t)_{t\geq 0}$. If we know
that $p_{t}^{n}$ converges to some $p_{t}$ in a sufficiently strong sense,
of course we obtain $P_{t}(x,dy)=p_{t}(x,y)dy.$ But in our framework $%
p_{t}^{n}$ does not converge: here, $p_{t}^{n}$ can even ``blow up'' as $%
n\rightarrow \infty $. However, if we may find a good equilibrium between
the blow-up and the speed of convergence, then we are able to conclude that $%
P_{t}(x,dy)=p_{t}(x,y)dy$ and $p_{t}$ has some regularity properties. This
is an interpolation type result.

Roughly speaking our main result is as follows. We assume that the speed of
convergence is controlled in the following sense: there exists some $a\in {%
\mathbb{N}}$ such that for every $q\in {\mathbb{N}}$ 
\begin{equation}
\left\Vert (L-L_{n})f\right\Vert _{q,\infty }\leq \varepsilon _{n}\left\Vert
f\right\Vert _{q+a,\infty }.  \label{i1}
\end{equation}%
Here $\left\Vert f\right\Vert _{q,\infty }$ is the norm in the standard
Sobolev space $W^{q,\infty }.$ In fact we will work with weighted Sobolev
spaces, and this is an important point. And also, we will assume a similar
hypothesis for the adjoint $(L-L_{n})^{\ast }$ (see Assumption \ref{A1A*1} 
for a precise statement).

Moreover we assume a ``propagation of regularity'' property: there exist $%
b\in {\mathbb{N}}$ and $\Lambda _{n}\geq 1$ such that for every $q\in {%
\mathbb{N}}$%
\begin{equation}
\left\Vert P_{t}^{n}f\right\Vert _{q,\infty }\leq \Lambda _{n}\left\Vert
f\right\Vert _{q+b,\infty }.  \label{i2}
\end{equation}%
Here also we will work with weighted Sobolev norms. And a similar hypothesis
is supposed to hold for the adjoint $P_{t}^{\ast ,n}$ (see Assumption \ref%
{A2A*2} 
for a precise statement).

Finally we assume the following regularity property: for every $t\in (0,1]$, 
$P_{t}^{n}(x,dy)=p_{t}^{n}(x,y)dy$ with $p_{t}^{n}\in C^{\infty }({\mathbb{R}%
}^{d}\times {\mathbb{R}}^{d})$ and for every $\kappa \geq 0$, $t\in (0,1]$, 
\begin{equation}
\left\vert \partial _{x}^{\alpha }\partial _{y}^{\beta
}p_{t}^{n}(x,y)\right\vert \leq \frac{C}{(\lambda _{n}t)^{\theta
_{0}(\left\vert \alpha \right\vert +\left\vert \beta \right\vert +\theta
_{1})}}\times \frac{(1+\left\vert x\right\vert ^{2})^{\pi (\kappa )}}{%
(1+\left\vert x-y\right\vert ^{2})^{\kappa }}.  \label{i3}
\end{equation}%
Here, $\alpha ,\beta $ are multi-indexes and $\partial _{x}^{\alpha
},\partial _{y}^{\beta }$ are the corresponding differential operators.
Moreover, $\pi (\kappa )$, $\theta _{0}$ and $\theta _{1}$ are suitable
parameters and $\lambda _{n}\rightarrow 0$ as $n\rightarrow \infty $ (we
refer to Assumption \ref{A3}).
In concrete examples (jump type stochastic differential equations) $\lambda
_{n}$ is related to the lower eigenvalue of the Malliavin covariance matrix
- essentially this is of order $\lambda _{n}^{\theta _{0}}.$ And in order to
handle the derivatives $\partial _{x}^{\alpha }$ and $\partial _{y}^{\beta }$
we need to make $\left\vert \alpha \right\vert $ respectively $\left\vert
\beta \right\vert $ integrations by parts (which involve $\lambda
_{n}^{\theta _{0}}).$ See also Assumption \ref{HH3}. 

By (\ref{i1})--(\ref{i3}), the rate of convergence is controlled by $%
\varepsilon _{n}\rightarrow 0$ and the blow-up of $p_{t}^{n}$ is controlled
by $\lambda _{n}^{-\theta _{0}}\rightarrow \infty $. So the regularity
property may be lost as $n\rightarrow \infty $. However, if there is a good
equilibrium between $\varepsilon _{n}\rightarrow 0$ and $\lambda
_{n}^{-\theta _{0}}\rightarrow \infty $ and $\Lambda _{n}\rightarrow \infty $
then the regularity is saved: we ask that for some $\delta >0$ 
\begin{equation}
\limsup_{n\to\infty}\frac{\varepsilon _{n}\Lambda _{n}}{\lambda _{n}^{\theta
_{0}(a+b+\delta )}}<\infty ,  \label{i4}
\end{equation}%
the parameters $a$, $b$ and $\theta _{0}$ being given in (\ref{i1}), (\ref%
{i2}) and (\ref{i3}) respectively. Then $P_{t}(x,dy)=p_{t}(x,y)dy$ with $%
p_{t}\in C^{\infty }({\mathbb{R}}^{d}\times {\mathbb{R}}^{d})$ and the
following upper bound holds: 
 for every $\varepsilon >0$, $\kappa \in {%
	\mathbb{N}}$ and $R>0$, one may find some constant $C, \pi(\kappa)>0$ such that for
very $(x,y) \in {\mathbb{R}}^d\times {\mathbb{R}}^{d}$ with $|x|<R$ and $t\in(0,1]$%
\begin{equation}
\left\vert \partial _{x}^{\alpha }\partial _{y}^{\beta
}p_{t}(x,y)\right\vert \leq \frac{C}{t^{\theta _{0}(1+\frac{a+b}{\delta }%
)(\left\vert \alpha \right\vert +\left\vert \beta \right\vert
+2d+\varepsilon )}}\times \frac{(1+|x|^2)^{\pi(\kappa)}}{(1+\left\vert
x-y\right\vert ^{2})^{\kappa }}.  \label{i3'}
\end{equation}%
This is the \textquotedblleft transfer of regularity\textquotedblright\ that
we mention in the title and which is stated in Theorem \ref{TransferBIS}.
The proof is based on a criterion of regularity for probability measures
given in \cite{[BC]}, which is close to interpolation spaces techniques.

The regularity criterion presented in this paper is tailored in order to
handle the following example (which will be treated in a forthcoming paper).
We consider the integro-differential operator 
\begin{equation}
Lf(x)=\<b(x),\nabla f(x)\>+\int_{E}\big(f(x+c(z,x))-f(x)-\<c(z,x),\nabla
f(x)\>\big)d\mu (z)  \label{i5}
\end{equation}%
where $\mu $ is an infinite measure on the normed space $(E,\vert \cdot
\vert _{E})$ such that $\int_{E}1\wedge \left\vert c(z,x)\right\vert
^{2}d\mu (z)<\infty .$ Moreover, for a sequence $\delta _{n}\downarrow 0 $,
we denote%
\begin{equation*}
A_{n}^{i,j}(x)=\int_{\{\left\vert z\right\vert _{E}\leq \delta
_{n}\}}c^{i}(z,x)c^{j}(z,x)d\mu (z)
\end{equation*}%
and we define%
\begin{equation}
\begin{array}{rl}
L_{n}f(x)= & \displaystyle\<b(x),\nabla f(x)\>+\int_{\{\left\vert
z\right\vert _{E}\geq \delta _{n}\}}(f(x+c(z,x))-f(x)-\<c(z,x),\nabla
f(x)\>)d\mu (z)\smallskip \\ 
& \displaystyle+\frac{1}{2}\mathrm{tr}(A_{n}(x)\nabla ^{2}f(x)).%
\end{array}
\label{i6}
\end{equation}%
By Taylor's formula, 
\begin{equation*}
\left\Vert Lf-L_{n}f\right\Vert _{\infty }\leq \left\Vert f\right\Vert
_{3,\infty }\varepsilon _{n}\quad \mbox{with}\quad \varepsilon
_{n}=\sup_{x}\int_{\{\left\vert z\right\vert _{E}\leq \delta
_{n}\}}\left\vert c(z,x)\right\vert ^{3}d\mu (z)
\end{equation*}%
(recall that $\|\cdot\|_{3,\infty}$ is the norm in the standard Sobolev
space $W^{3,\infty}$). Under the uniform ellipticity assumption $%
A_{n}(x)\geq \lambda _{n}$ for every $x\in {\mathbb{R}}^{d},$ the semigroup $%
(P^n_t)_{t\geq 0}$ associated to $L_{n} $ has the regularity property (\ref%
{i3}) with $\theta _{0}$ depending on the measure $\mu .$ The speed of
convergence in (\ref{i1}), with $a=3,$ is controlled by $\varepsilon
_{n}\downarrow 0.$ So, if (\ref{i4}) holds, then we obtain the regularity of 
$P_{t}$ and the short time estimates (\ref{i3'}).

The semigroup $(P_t)_{t\geq 0}$ associated to $L$ corresponds to stochastic
equations driven by the Poisson point measure $N_{\mu }(dt,dz)$ with
intensity measure $\mu $, so the problem of the regularity of $P_{t}$ has
been extensively discussed in the probabilistic literature. A first approach
initiated by Bismut \cite{[Bi]}, L\'eandre \cite{[L]} and Bichteler,
Gravereaux and Jacod \cite{[BGJ]} (see also the recent monograph of Bouleau
and Denis \cite{[BD]} and the bibliography therein), is done under the
hypothesis that $E={\mathbb{R}}^{m}$ and $\mu (dz)=h(z)dz$ with $h\in
C^{\infty }({\mathbb{R}}^{m}).$ Then one constructs a Malliavin type
calculus based on the amplitude of the jumps of the Poisson point measure $%
N_{\mu }$ and employs this calculus in order to study the regularity of $%
P_{t}.$ A second approach initiated by Carlen and Pardoux \cite{[CP]} (see
also Bally and Cl\'{e}ment \cite{[BCl]}) follows the ideas in Malliavin
calculus based on the exponential density of the jump times in order to
study the same problem. Finally a third approach is due to Picard \cite%
{[P1],[P2]}, but see also Ishikawa and Kunita \cite{[IK]}, the contributions
of Kunita \cite{[K1],[K2]} and the recent monograph by Ishikawa \cite{[I]}
for many references and developments in this direction. Picard constructs a
Malliavin type calculus based on finite differences (instead of standard
Malliavin derivatives) and obtains the regularity of $P_{t}$ for a general
class of intensity measures $\mu $ including purely atomic measures (in
contrast with $\mu (dz)=h(z)dz)$. We stress that all the above approaches
work under different non degeneracy hypotheses, each of them corresponding
to the specific noise that is used in the calculus. So in some sense we have
not a single problem but three different classes of problems. The common
feature is that the strategy in order to solve the problem follows the ideas
from Malliavin calculus based on some noise contained in $N_{\mu }.$ Our
approach is completely different because, as described above, we use the
regularization effect of $\mathrm{tr}(A_{n}(x)\nabla ^{2}).$ This
regularization effect may be exploited either by using the standard
Malliavin calculus based on the Brownian motion or using some analytical
arguments. The approach that we propose in \cite{[BCW]} is probabilistic, so
employs the standard Malliavin calculus. But anyway, as mentioned above, the
regularization effect vanishes as $n\rightarrow \infty $ and a supplementary
argument based on the equilibrium given in (\ref{i4}) is used. We precise
that the non degeneracy condition $A_{n}(x)\geq \lambda _{n}>0$ is of the
same nature as the one employed by J. Picard so the problem we solve is in
the same class.

The idea of replacing \textquotedblleft small jumps\textquotedblright\ (the
ones in $\{\left\vert z\right\vert _{E}\leq \varepsilon _{n}\}$ here) by a
Brownian part (that is $\mathrm{tr}(A_{n}(x)\nabla ^{2})$ in $L_{n})$ is not
new -- it has been introduced by Asmussen and Rosinski in \cite{[AR]} and
has been extensively employed in papers concerned with simulation problems:
since there is a huge amount of small jumps, they are difficult to simulate
and then one approximates them by the Brownian part corresponding to $%
\mathrm{tr}(A_{n}(x)\nabla ^{2}).$ See for example \cite{[A],[BK],[CD]} and
many others. However, at our knowledge, this idea has not been yet used in
order to study the regularity of $P_{t}.$

The paper is organized as follows. In Section \ref{sect:NotRes} we give the
notation and the main results mentioned above and in Section \ref%
{sect:proofs} we give the proof of these results. Section \ref{sect:reg} is
devoted to some preliminary results about regularity. Namely, in Section \ref%
{sect:3.1} we recall and develop some results concerning regularity of
probability measures, based on interpolation type arguments, coming from 
\cite{[BC]}. These are the main instruments used in the paper. In Section %
\ref{sect:3.2} we prove a regularity result which is a key point in our
approach. In fact, it allows to handle the multiple integrals coming from
the application of a Lindeberg method for the decomposition of $P_t-P^n_t$.
Finally, in Appendix \ref{app:weights} and \ref{app:semi} we prove some
technical results used in the paper.

\section{Notation and main results}

\label{sect:NotRes}

\subsection{Notation}

\label{sect:notation}

For a multi-index $\alpha =(\alpha _{1},...,\alpha _{m})\in \{1,...,d\}^{m}$
we denote $\left\vert \alpha \right\vert =m$ (the length of the multi-index)
and $\,\partial ^{\alpha }$ is the derivative corresponding to $\alpha ,$
that is $\partial ^{\alpha _{m}}\cdots \partial ^{\alpha _{1}}$, with $%
\partial ^{\alpha _{i}}=\partial _{x_{\alpha _{i}}}$. For $f\in C^{\infty }({%
\mathbb{R}}^{d}\times {\mathbb{R}}^{d})$, $(x,y)\in {\mathbb{R}}^{d}\times {%
\mathbb{R}}^{d}$ and two multi-indexes $\alpha $ and $\beta $, we denote by $%
\partial _{x}^{\alpha }$ the derivative with respect to $x$ and by $\partial
_{y}^{\beta }$ the derivative with respect to $y$. Moreover, for $f\in
C^{\infty }({\mathbb{R}}^{d})$ and $q\in {\mathbb{N}}$ we denote%
\begin{equation}
\left\vert f\right\vert _{q}(x)=\sum_{0\leq \left\vert \alpha \right\vert
\leq q}\left\vert \partial ^{\alpha }f(x)\right\vert .  \label{NOT1}
\end{equation}%
If $f$ is not a scalar function, that is, $f=(f^{i})_{i=1,\ldots ,d}$ or $%
f=(f^{i,j})_{i,j=1,\ldots ,d}$, we denote $\left\vert f\right\vert
_{q}=\sum_{i=1}^{d}\left\vert f^{i}\right\vert _{q}$ respectively $%
\left\vert f\right\vert _{q}=\sum_{i,j=1}^{d}\left\vert f^{i,j}\right\vert
_{q}.$

We will work with the weights 
\begin{equation}
\psi _{\kappa }(x)=(1+\left\vert x\right\vert ^{2})^{\kappa },\quad \kappa
\in {\mathbb{Z}} .  \label{NOT2}
\end{equation}%
The following properties hold:

\begin{itemize}
\item for every $\kappa\geq \kappa^{\prime }\geq 0$, 
\begin{equation}
\psi _{\kappa }(x) \leq \psi _{\kappa ^{\prime }}(x);  \label{NOT3a}
\end{equation}

\item for every $\kappa\geq 0$, there exists $C_{\kappa }>0$ such that 
\begin{equation}
\psi _{\kappa }(x)\leq C_{\kappa }\psi _{\kappa }(y)\psi _{\kappa }(x-y);
\label{NOT3b}
\end{equation}

\item for every $\kappa \geq 0$, there exists $C_{\kappa }>0$ such that for
every $\phi \in C_{b}^{\infty }({\mathbb{R}}^{d})$, 
\begin{equation}
\psi _{\kappa }(\phi (x))\leq C_{\kappa }\psi _{\kappa }(\phi
(0))(1+\left\Vert \nabla \phi \right\Vert _{\infty }^{2})^{\kappa }\psi
_{\kappa }(x);  \label{NOT3d}
\end{equation}

\item for every $q\in {\mathbb{N}}$ there exist $\overline{C}_{q}, 
\underline{C}_{q}>0$ such that for every $\kappa \in {\mathbb{R}}$ and $f\in
C^{\infty }({\mathbb{R}}^{d})$, 
\begin{equation}
\underline{C}_{q}\psi _{\kappa }\left\vert f\right\vert _{q}(x)\leq
\left\vert \psi _{\kappa }f\right\vert _{q}(x)\leq \overline{C}_{q}\psi
_{\kappa }\left\vert f\right\vert _{q}(x).  \label{NOT3c}
\end{equation}
\end{itemize}

%
%
Note that (\ref{NOT3a})--(\ref{NOT3d}) are immediate, whereas (\ref{NOT3c})
is proved in Appendix \ref{app:weights} (see Lemma \ref{Psy1}).

For $q\in {\mathbb{N}}$, $\kappa \in {\mathbb{R}}$ and $p\in (1,\infty ]$
(we stress that we include the case $p=+\infty $), we set $\Vert \cdot \Vert
_{p}$ the usual norm in $L^{p}({\mathbb{R}}^{d})$ and 
\begin{equation}
\left\Vert f\right\Vert _{q,\kappa ,p}=\big\| \vert \psi _{\kappa }f\vert
_{q}\big\| _{p}.  \label{NOT4}
\end{equation}%
We denote $W^{q,\kappa ,p}$ to be the closure of $C^{\infty }({\mathbb{R}}%
^{d})$ with respect to the above norm. If $\kappa =0$ we just denote $%
\left\Vert f\right\Vert _{q,p}=\left\Vert f\right\Vert _{q,0,p}$ and $%
W^{q,p}=W^{q,0,p}$ (which is the usual Sobolev space). So, we are working
with weighted Sobolev spaces. 
The weighted Sobolev spaces $W^{q,\kappa ,p}$ are the natural framework in
the paper \cite{[BC]} where the ``balance argument'' is obtained. There (see Theorem A.2 in \cite{[BC]}) we
have used a crucial result of Petrushev and Xu \cite{[PX]} concerning the
construction of kernels with polynomial decay at infinity. Then the weights 
$\psi _{\kappa }$ appear in a natural way in order to capture the behaviour 
of the kernel at infinity.

The following properties hold:

\begin{itemize}
\item for every $q\in {\mathbb{N}}$ there exists $\overline{C}_{q}\geq 
\underline{C}_{q}>0$ such that for every $\kappa \in {\mathbb{R}}$, $p>1$
and $f\in W^{q,\kappa,p}({\mathbb{R}}^{d})$, 
\begin{equation}
\underline{C}_{q}\Vert \psi _{\kappa }|f|_{q}\Vert _{p}\leq \Vert f\Vert
_{q,\kappa ,p}\leq \overline{C}_{q}\Vert \psi _{\kappa }|f|_{q}\Vert _{p};
\label{NOT4a}
\end{equation}

\item for every $q\in {\mathbb{N}}$ and $p>1$ there exists $C_{q,p}>0$ such
that for every $\kappa \in {\mathbb{R}}$ and $f\in W^{q,k,p}({\mathbb{R}}%
^{d})$, 
\begin{equation}
\left\Vert f\right\Vert _{q,\kappa ,p}\leq C_{q,p}\left\Vert f\right\Vert
_{q,\kappa +d,\infty }  \label{NOT5a}
\end{equation}%
and if $p>d$, 
\begin{equation}
\left\Vert f\right\Vert _{q,\kappa ,\infty }\leq C_{q,p}\left\Vert
f\right\Vert _{q+1,\kappa ,p};  \label{NOT5b}
\end{equation}

\item for $\kappa ,\kappa ^{\prime }\in {\mathbb{R}}$, $q,q^{\prime }\in {%
\mathbb{N}}$, $p\in (1,\infty ]$ and $U:C^{\infty }({\mathbb{R}}%
^{d})\rightarrow C^{\infty }({\mathbb{R}}^{d})$, the following two
assertions are equivalent: there exists a constant $C_{\ast }\geq 1$ such
that for every $f$, 
\begin{equation}
\left\Vert Uf\right\Vert _{q,\kappa ,\infty }\leq C_{\ast }\left\Vert
f\right\Vert _{q^{\prime },\kappa ^{\prime },p}  \label{NOT6a}
\end{equation}%
and there exists a constant $C^{\ast }\geq 1$ such that for every $f$, 
\begin{equation}
\Big\Vert\psi _{\kappa }U\Big(\frac{1}{\psi _{\kappa ^{\prime }}}f\Big)%
\Big\Vert_{q,\infty }\leq C^{\ast }\left\Vert f\right\Vert _{q^{\prime },p} .
\label{NOT6b}
\end{equation}
\end{itemize}

%
Notice that (\ref{NOT4a}) is a consequence of (\ref{NOT3c}). The inequality (%
\ref{NOT5a}) is an immediate consequence of (\ref{NOT3c}) and of the fact
that $\psi _{-d}\in L^{p}({\mathbb{R}}^{d})$ for every $p\geq 1$. And the
inequality (\ref{NOT5b}) is a consequence of Morrey's inequality (Corollary
IX.13 in \cite{Morrey}), whose use gives $\left\Vert f\right\Vert
_{0,0,\infty }\leq \left\Vert f\right\Vert _{1,0,p}$, and of (\ref{NOT3c}).
In order to prove the equivalence between (\ref{NOT6a}) and (\ref{NOT6b}),
one takes $g=\psi _{\kappa ^{\prime }}f$ (respectively $g=\frac{1}{\psi
_{\kappa ^{\prime }}}f)$ and uses (\ref{NOT3c}) as well.

\subsection{Main results}

\label{sect:results}

We consider a Markov semigroup $(P_t)_{t\geq 0}$ with infinitesimal operator 
$L$ and a sequence $(P^n_t)_{t\geq 0}$, $n\in {\mathbb{N}}$, of Markov
semigroups with infinitesimal operator $L_{n}.$ We suppose that ${\mathcal{S}%
({\mathbb{R}}^{d})}$ is included in the domain of $(P_t)_{t\geq 0}$, $%
(P^n_t)_{t\geq 0}$, $L$ and of $L_{n}$ and we suppose that for $f\in {%
\mathcal{S}({\mathbb{R}}^{d})}$ 
we have $P_tf, P^n_tf, Lf, L_{n}f \in {\mathcal{S}({\mathbb{R}}^{d})}$.

We denote $\Delta _{n}=L-L_{n}.$ Moreover, we denote by $P_{t}^{\ast ,n}$
the formal adjoint of $P_{t}^{n}$ and by $\Delta _{n}^{\ast }$ the formal
adjoint of $\Delta _{n}$ that is%
\begin{equation}
\left\langle P_{t}^{\ast ,n}f,g\right\rangle =\left\langle
f,P_{t}^{n}g\right\rangle \quad \mbox{and}\quad \left\langle \Delta
_{n}^{\ast }f,g\right\rangle =\left\langle f,\Delta _{n}g\right\rangle ,
\label{TR1}
\end{equation}%
$\left\langle \cdot ,\cdot \right\rangle $ being the scalar product in $%
L^{2}({\mathbb{R}}^{d},dx).$

We present now our hypotheses. The first one concerns the speed of
convergence of $L_{n}\rightarrow L.$

\begin{assumption}
\label{A1A*1} Let $a\in {\mathbb{N}}$, and let $(\varepsilon _{n})_{n\in {%
\mathbb{N}}}$ be a decreasing sequence such that $\lim_{n\to\infty}%
\varepsilon _{n}=0.$We assume that for every $q\in {\mathbb{N}},\kappa \geq 0
$ and $p>1$ there exists $C>0$ such that for every $n\in{\mathbb{N}}$ and $%
f\in\mathcal{S}({\mathbb{R}}^d)$, 
\begin{align}
(A_{1})& \qquad \left\Vert \Delta _{n}f\right\Vert _{q,-\kappa ,\infty }\leq
C\varepsilon _{n}\left\Vert f\right\Vert _{q+a,-\kappa ,\infty },
\label{TR3} \\
(A_{1}^{\ast })& \qquad \left\Vert \Delta _{n}^{\ast }f\right\Vert
_{q,\kappa ,p}\leq C\varepsilon _{n}\left\Vert f\right\Vert _{q+a,\kappa ,p}.
\label{TR3'}
\end{align}
\end{assumption}

Our second hypothesis concerns the ``propagation of regularity'' for the
semigroups $(P^n_t)_{t\geq 0}$.

\begin{assumption}
\label{A2A*2} Let $\Lambda _{n}\geq 1,n\in {\mathbb{N}},$ be an increasing
sequence such that $\Lambda _{n+1}\leq \gamma \Lambda _{n}$ for some $\gamma
\geq 1.$ For every $q\in {\mathbb{N}}$ and $\kappa \geq 0,p>1$, there exist $%
C>0$ and $b\in {\mathbb{N}}$, such that for every $n\in {\mathbb{N}}$ and $%
f\in \mathcal{S}({\mathbb{R}}^{d})$%
\begin{align}
(A_{2})& \qquad \sup_{s\leq t}\left\Vert P_{s}^{n}f\right\Vert _{q,-\kappa
,\infty }\leq C\Lambda _{n}\left\Vert f\right\Vert _{q+b,-\kappa ,\infty },
\label{TR2} \\
(A_{2}^{\ast })& \qquad \sup_{s\leq t}\left\Vert P_{s}^{\ast ,n}f\right\Vert
_{q,\kappa ,p}\leq C\Lambda _{n}\left\Vert f\right\Vert _{q+b,\kappa ,p}.
\label{TR2'}
\end{align}
\end{assumption}

The hypothesis $(A_{2}^{\ast })$ is rather difficult to verify so, in
Appendix \ref{app:semi}, we give some sufficient conditions in order to
check it (see Proposition \ref{A2}).

Our third hypothesis concerns the ``regularization effect'' of the
semi-group $(P^n_t)_{t\geq 0}$.

\begin{assumption}
\label{A3} We assume that 
\begin{equation}
P_{t}^{n}f(x)=\int_{{\mathbb{R}}^{d}}p_{t}^{n}(x,y)f(y)dy  \label{TR4}
\end{equation}%
with $p_{t}^{n}\in C^{\infty }({\mathbb{R}}^{d}\times {\mathbb{R}}^{d})$.
Moreover, we assume there exist $\theta _{0}>0$ and a sequence $\lambda _{n}$%
, $n\in {\mathbb{N}}$, with, as $n\to\infty$, 
\begin{equation}
\lambda _{n}\downarrow 0,\quad \lambda _{n}\leq \gamma \lambda _{n+1},
\label{TRa1-new}
\end{equation}%
for some $\gamma \geq 1,$ such that the following property holds: for every $%
\kappa \geq 0,q\in {\mathbb{N}}$ there exist $\pi (q,\kappa ),$ increasing
in $q$ and in $\kappa ,$ a constant $\theta _{1}\geq 0,$ and a constant $C>0$
such that for every $n\in {\mathbb{N}}$, $t\in (0,1]$, for every
multi-indexes $\alpha $ and $\beta $ with $\left\vert \alpha \right\vert
+\left\vert \beta \right\vert \leq q$ and $(x,y)\in {\mathbb{R}}^{d}\times {%
\mathbb{R}}^{d}$ 
\begin{equation}
(A_{3})\quad \left\vert \partial _{x}^{\alpha }\partial _{y}^{\beta
}p_{t}^{n}(x,y)\right\vert \leq C\frac{1}{(\lambda _{n}t)^{\theta
_{0}(q+\theta _{1})}}\times \frac{\psi _{\pi (q,\kappa )}(x)}{\psi _{\kappa
}(x-y)}.  \label{TR5}
\end{equation}
\end{assumption}

Note that in (\ref{TR5}) we are quantifying the possible blow-up of $%
|\partial _{x}^{\alpha }\partial _{y}^{\beta }p_{t}^{n}(x,y)|$ as $%
n\rightarrow \infty $.

We also assume the following statements hold for the semigroup $(P_t)_{t\geq
0}$.

\begin{assumption}
\label{A5} For every $\kappa \geq 0,k\in {\mathbb{N}}$ there exists $C\geq 1$
such that 
\begin{equation}
(A_{4})\quad \left\Vert P_{t}f\right\Vert _{k,-\kappa ,\infty }\leq
C\left\Vert f\right\Vert _{k,-\kappa ,\infty }.  \label{R7}
\end{equation}
\end{assumption}

\begin{assumption}
\label{A6} For every $\kappa \geq 0,k\in {\mathbb{N}}$ there exists $C\geq 1,%
\overline{\kappa }\geq \kappa $ such that 
\begin{equation}
(A_{5})\quad P_{t}\psi _{\kappa }(x)\leq C\psi _{\overline{\kappa }}(x).
\label{R7'''}
\end{equation}
\end{assumption}

Our main result is the following:

\begin{theorem}
\label{TransferBIS} Suppose that Assumption \ref{A1A*1}, \ref{A2A*2} , \ref%
{A3}, \ref{A5} and \ref{A6} 
hold. Suppose also that for some $\delta >0$ 
\begin{equation*}
\limsup_{n\to\infty}\frac{\varepsilon _{n}\Lambda _{n}}{\lambda _{n}^{\theta
_{0}(a+b+\delta )}}<\infty.
\end{equation*}%
Then $P_{t}(x,y)=p_{t}(x,y)$ with $p_{t}\in C^{\infty }({\mathbb{R}}%
^{d}\times {\mathbb{R}}^{d}).$ Moreover, for every $\kappa \in {\mathbb{N}}%
,R\in {\mathbb{N}},\varepsilon >0$ and every multi-indexes $\alpha $ and $%
\beta $ there exists some constants $C=C(\delta,R,\kappa ,\varepsilon ,\alpha
,\beta )$ and $t_0=t_0(\delta,R,\kappa ,\varepsilon ,\alpha
,\beta )\in(0,1]$ such that for every $t\in(0,t_0)$, $x\in\R^d$ with $|x|<R$ and $y\in {\mathbb{R}}%
^{d}$, 
\begin{equation}
\left\vert \partial _{x}^{\alpha }\partial _{y}^{\beta
}p_{t}(x,y)\right\vert \leq C\times t^{-\theta _{0}(1+\frac{a+b}{\delta }%
)(\left\vert \alpha \right\vert +\left\vert \beta \right\vert
+2d+\varepsilon )}\times \frac{1}{\psi _{\kappa }(x-y)}  \label{TR6e}
\end{equation}%
with $\theta _{0}$ from (\ref{TR5}).
\end{theorem}

\section{Regularity results}

\label{sect:reg}

This section is devoted to some preliminary results allowing us to prove the
statements resumed in Section \ref{sect:results}: in Section \ref{sect:3.1}
we give an abstract regularity criterion, while in Section \ref{sect:3.2} we
prove a regularity result for iterated integrals, that will be useful to
handle a Lindeberg type decomposition of $P_t-P^n_t$.

\subsection{A regularity criterion based on interpolation}

\label{sect:3.1}

Let us first recall some results obtained in \cite{[BC]} concerning the
regularity of a measure $\mu $ on ${\mathbb{R}}^{d}$ (with the Borel $\sigma 
$-field). For two signed finite measures $\mu ,\nu $ and for $k\in {\mathbb{N%
}}$ we define the distance%
\begin{equation}
d_{k}(\mu ,\nu )=\sup \Big\{\Big\vert\int fd\mu -\int fd\nu \Big\vert:\Vert
f\Vert _{k,\infty }\leq 1\Big\}.  \label{reg2}
\end{equation}%
If $\mu $ and $\nu $ are probability measures, $d_{0}$ is the total
variation distance and $d_{1}$ is the Fortet Mourier distance. In this paper
we will work with an arbitrary $k\in {\mathbb{N}}$. Notice also that $%
d_{k}(\mu ,\nu )=\left\Vert \mu -\nu \right\Vert _{W_{\ast }^{k,\infty }}$
where $W_{\ast }^{k,\infty }$ is the dual of $W^{k,\infty }.$

We fix now $k,q,h\in {\mathbb{N}}$, with $h\geq 1$, and $p>1$. Hereafter, we
denote by $p_{\ast }=p/(p-1)$ the conjugate of $p.$ Then, for a signed
finite measure $\mu $ and for a sequence of absolutely continuous signed
finite measures $\mu _{n}(dx)=f_{n}(x)dx$ with $f_{n}\in C^{2h+q}({\mathbb{R}%
}^{d}),$ we define%
\begin{equation}
\pi _{k,q,h,p}(\mu ,(\mu _{n})_{n})=\sum_{n=0}^{\infty }2^{n(k+q+d/p_{\ast
})}d_{k}(\mu ,\mu _{n})+\sum_{n=0}^{\infty }\frac{1}{2^{2nh}}\left\Vert
f_{n}\right\Vert _{2h+q,2h,p}.  \label{reg3}
\end{equation}

The following result is the key point in our approach:

\begin{lemma}
\label{lemma-inter} Let $k,q,h\in {\mathbb{N}}$ with $h\geq 1$, and $p>1$ be
given. There exists a constant $C_*$ (depending on $k,q,h$ and $p$ only)
such that the following holds. Let $\mu $ be a finite measure for which one
may find a sequence $\mu _{n}(dx)=f_{n}(x)dx$, $n\in {\mathbb{N}}$ such that 
$\pi _{k,q,h,p}(\mu ,(\mu _{n})_{n})<\infty .$ Then $\mu (dx)=f(x)dx$ with $%
f\in W^{q,p}$ and moreover 
\begin{equation}
\left\Vert f\right\Vert _{q,p}\leq C_{\ast }\times \pi _{k,q,h,p}(\mu ,(\mu
_{n})_{n}).  \label{reg4}
\end{equation}
\end{lemma}

The proof of Lemma \ref{lemma-inter} is given in \cite{[BC]}, being a
particular case (take $\mathbf{e}=\mathbf{e}_{p}$) of Proposition A.1 in
Appendix A. We give a first simple consequence:

\begin{lemma} 
	\label{reg} Let $p_{t}\in C^{\infty }({\mathbb{R}}^{d}),t>0,$ be a family of
	probability densities such that, for every $\kappa \in {\mathbb{N}}$,  $\int \psi _{\kappa }(x)p_{t}(x)dx\leq
	m_{\kappa }<\infty $ . We assume that for
	some $\theta _{0}>0$ and $\theta _{1}>0$ the following holds. For every $%
	q,\kappa \in {\mathbb{N}}$ and $p\geq 1$ there exists a constant $%
	C=C(q,\kappa ,p, \theta_1)\geq 1$ such that 
	\begin{equation}
	\left\Vert \psi _{\kappa }p_{t}\right\Vert _{q,p}\leq Ct^{-\theta
		_{0}(q+\theta _{1})} .  \label{a1}
	\end{equation}%
	Then, for every $\delta >0$ there exists a constant $C=C(q,\kappa ,p,\delta, \theta_1
	) $ such that 
	\begin{equation}
	\left\Vert \psi _{\kappa }p_{t}\right\Vert _{q,p}\leq Ct^{-\theta _{0}(q+%
		\frac{d}{p_{\ast }}+\delta)}  \label{a2}
	\end{equation}%
	for every $t<2^ {-4/\theta_0}$ (so does not matter the value of $\theta _{1},$ one may morally  replace it
	by $\frac{d}{p_{\ast }}$ in the power of $t$; however, $\theta_1$ appears in the constant $C$).
\end{lemma}

\textbf{Proof.} We take $n_{\ast }\in {\mathbb{N}}$ and we define $f_{n}=0$
for $n\leq n_{\ast }$ and $f_{n}=\psi_{\kappa}p_{t}$ for $n>n_{\ast }.$
Notice that $d_{0}(\psi _{\kappa }p_{t},0)=m_{\kappa }.$ Then (\ref{reg4})
with $k=0$ gives 
\begin{eqnarray*}
	\left\Vert \psi _{\kappa }p_{t}\right\Vert _{q,p} &\leq &C_{\ast }\Big(%
	m_{\kappa}\sum_{n=0}^{n_{\ast }}2^{n(q+\frac{d}{p_{\ast }})}+\left\Vert \psi
	_{\kappa }p_{t}\right\Vert _{2h+q,2h,p}\sum_{n=n_{\ast }+1}^{\infty }\frac{1%
	}{2^{2nh}}\Big) \\
	&\leq &C_{\ast }\Big(m_{\kappa}2^{n_{\ast }(q+\frac{d}{p_{\ast }}%
		)}+\left\Vert \psi _{\kappa }p_{t}\right\Vert _{2h+q,2h,p}\frac{1}{%
		2^{2n_{\ast }h}}\Big).
\end{eqnarray*}%
Since $\| \psi _{\kappa }p_{t}\| _{2h+q,2h,p}=\| \psi _{2h+\kappa }p_{t}\|
_{2h+q,p}$, we use (\ref{a1}) (with $q$ replaced by $2h+q$ and $\kappa$ replaced by $2h+\kappa$) and we obtain
\begin{equation}\label{new}
\left\Vert \psi _{\kappa }p_{t}\right\Vert _{q,p} \leq 
C_{\ast }\Big(m_{\kappa}2^{n_{\ast }(q+\frac{d}{p_{\ast }}%
	)}+C_ht^{-\theta
	_{0}(2h+q+\theta _{1})} \frac{1}{%
	2^{2n_{\ast }h}}\Big).
\end{equation}%
We write here $C_h$ to stress the possible dependence on $h$ of the constant $C$ in (\ref{a1}). We optimize over $n_{\ast }
$: we look for $n_\ast$ such that 
$$
m_{\kappa}2^{n_{\ast }(q+\frac{d}{p_{\ast }}%
	)}=C_ht^{-\theta
	_{0}(2h+q+\theta _{1})} \frac{1}{%
	2^{2n_{\ast }h}}.
$$
Straightforward computations give $n_\ast=\lfloor n_\ast(h)\rfloor$ where 
$$
n_\ast(h)=\frac{\theta_0(2h+q+\theta_1)}{2h+q+d/p_\ast}\, \log_2t^{-1}- \frac{1}{2h+q+d/p_\ast}\,\log_2m_\kappa+\frac{1}{2h+q+d/p_\ast}\log_2C_h.
$$
In order to successfully insert such $n_\ast$ in (\ref{new}), we need that $n_\ast\geq 1$. Notice in fact that if $n_\ast=0$ then (\ref{new}) would give 
$
\left\Vert \psi _{\kappa }p_{t}\right\Vert _{q,p} \leq 
C_{\ast }(m_{\kappa}+C_ht^{-\theta
	_{0}(2h+q+\theta _{1})} )
$
and this is not the kind of estimate we are looking for. 

We take $h$ sufficiently large in order to have $2h\geq \log_2 m_\kappa$ and $\eta_h:=(q+\frac{d}{p_{\ast }})/2h\leq 1$. Then we write
$$
n_\ast(h)\geq \frac{\theta_0}{1+\eta_h}\, \log_2t^{-1}- \frac{1}{2h}\,\log_2m_\kappa\geq \frac{\theta_0}{2}\, \log_2t^{-1}-1.
$$
It follows that the restriction $n_\ast(h)\geq 1$ amounts to $\theta_0 \log_2t^{-1}\geq 4$ which is exactly the restriction we have assumed for $t$.

We replace $n_\ast(h)$ in (\ref{new}) and we obtain
\begin{eqnarray*}
	\left\Vert \psi _{\kappa }p_{t}\right\Vert _{q,p} 
	&\leq &2C_{\ast }\times m_{\kappa }^{\frac{1}{1+\eta_h}}\times
	C_h^\frac{\eta_h}{1+\eta_h} t^{-\theta _{0}(2h+q+\theta _{1})\frac{\eta_h}{1+\eta_h}}.
\end{eqnarray*}%
Since $(2h+q+\theta _{1})\frac{\rho
	_{h}}{1+\eta_h}\downarrow q+\frac{d}{p_{\ast }}$ as $h\to \infty$, $h$ can be chosen large enough in order that (\ref{a2}) holds (since $h$ depends on $\theta_1$ so does $C_h$ and consequently the constant in our estimates). $\square $

\medskip

We will also use the following consequence of Lemma \ref{lemma-inter} (the
proof is given in \cite{[B]} and we do not repeat it here):

\begin{lemma}
\label{REG} Let $k,q,h\in {\mathbb{N}}$, with $h\geq 1$, and $p>1$ be given
and set 
\begin{equation}
\rho _{h}:=\frac{k+q+d/p_{\ast }}{2h}.  \label{reg5}
\end{equation}%
We consider an increasing sequence $\theta (n)\geq 1,n\in {\mathbb{N}}$,
such that $\lim_{n\to\infty}\theta (n)=\infty $\ and $\theta (n+1)\leq
\Theta \times \theta (n)$ for some constant $\Theta \geq 1.$ Suppose that we
may find a sequence of functions $f_{n}\in C^{2h+q}({\mathbb{R}}^{d}),n\in {%
\mathbb{N}}$, such that 
\begin{equation}
\left\Vert f_{n}\right\Vert _{2h+q,2h,p}\leq \theta (n)  \label{reg9}
\end{equation}%
and, with $\mu _{n}(dx)=f_{n}(x)dx,$ 
\begin{equation}
\limsup_{n\to\infty}d_{k}(\mu ,\mu _{n})\times \theta ^{\rho
_{h}+\varepsilon }(n)<\infty  \label{reg10}
\end{equation}%
for some $\varepsilon >0.$ Then $\mu (dx)=f(x)dx$ with $f\in W^{q,p}.$

Moreover, for $\delta ,\varepsilon >0$ and $n_{\ast }\in {\mathbb{N}}$, let 
\begin{align}
A(\delta )& =\left\vert \mu \right\vert ({\mathbb{R}}^{d})\times 2^{l(\delta
)(1+\delta )(q+k+d/p_{\ast })}\quad 
\mbox{with $l(\delta )=\min
\{l:2^{l\times \frac{\delta }{1+\delta }}\geq l\}$},  \label{reg12'} \\
B(\varepsilon )& =\sum_{l=1}^{\infty }\frac{l^{2(q+k+d/p_{\ast }+\varepsilon
)}}{2^{2\varepsilon l}},  \label{reg12''} \\
C_{h,n_{\ast }}(\varepsilon )& =\sup_{n\geq n_{\ast }}d_{k}(\mu ,\mu
_{n})\times \theta ^{\rho _{h}+\varepsilon }(n).  \label{reg11}
\end{align}%
Then, for every $\delta >0$ 
\begin{equation}
\left\Vert f\right\Vert _{q,p}\leq C_{\ast }(\Theta +A(\delta )\theta
(n_{\ast })^{\rho _{h}(1+\delta )}+B(\varepsilon )C_{h,n_{\ast
}}(\varepsilon )),  \label{reg12}
\end{equation}%
$C_{\ast }$ being the constant in (\ref{reg4}) and $\rho _{h}$ being given
in (\ref{reg5}).
\end{lemma}

\subsection{A regularity lemma}

\label{sect:3.2}

We give here a regularization result in the following abstract framework. We
consider a sequence of operators $U_{j}:{\mathcal{S}({\mathbb{R}}^{d})}%
\rightarrow {\mathcal{S}({\mathbb{R}}^{d})}$, 
$j\in {\mathbb{N}}$, and we denote by $U_{j}^{\ast }$ the formal adjoint
defined by $\langle U_{j}^{\ast }f,g\rangle =\langle f,U_{j}g\rangle $ with
the scalar product in $L^{2}({\mathbb{R}}^{d})$.

\begin{assumption}
\label{H1H*1} Let $a\in {\mathbb{N}}$ be fixed. We assume that for every $%
q\in {\mathbb{N}},\kappa \geq 0$ and $p\in \lbrack 1,\infty )$ there exist
constants $C_{q,\kappa ,p}(U)$ and $C_{q,\kappa ,\infty }(U)$ such that for
every $j$ and $f$, 
\begin{align}
(H_{1})& \qquad \left\Vert U_{j}f\right\Vert _{q,-\kappa ,\infty }\leq
C_{q,\kappa ,\infty }(U)\left\Vert f\right\Vert _{q+a,-\kappa ,\infty },
\label{h1} \\
(H_{1}^{\ast })& \qquad \left\Vert U_{j}^{\ast }f\right\Vert _{q,\kappa
,p}\leq C_{q,\kappa ,p}(U)\left\Vert f\right\Vert _{q+a,\kappa ,p}.
\label{h1'}
\end{align}%
We assume that $C_{q,\kappa ,p}(U)$, $p\in \lbrack 1,\infty ]$, is non
decreasing with respect to $q$ and $\kappa $.
\end{assumption}

We also consider a semigroup $(S_t)_{t\geq 0}$ of the form 
\begin{equation*}
S_{t}(x,dy)=s_{t}(x,y)dy\quad \mbox{with}\quad s_{t}\in {\mathcal{S}}({%
\mathbb{R}}^{d}\times {\mathbb{R}}^{d}).
\end{equation*}

We define the formal adjoint operator 
\begin{equation*}
S_{t}^{\ast }f(y)=\int_{{\mathbb{R}}^{d}}s_{t}(x,y)f(x)dx,\quad t>0.
\end{equation*}

\begin{assumption}
\label{H2H*2} If $f\in {\mathcal{S}({\mathbb{R}}^{d})}$ 
then $S_{t}f\in {\mathcal{S}({\mathbb{R}}^{d})}$. 
Moreover, there exists $b\in {\mathbb{N}}$ such that for every $q\in {%
\mathbb{N}},\kappa \geq 0$ and $p\in \lbrack 1,\infty )$ there exist
constants $C_{q,\kappa ,p}(S)$ such that for every $t>0$, 
\begin{align}
(H_{2})& \qquad \left\Vert S_{t}f\right\Vert _{q,-\kappa ,\infty }\leq
C_{q,\kappa ,\infty }(S)\left\Vert f\right\Vert _{q+b,-\kappa ,\infty },
\label{h2} \\
(H_{2}^{\ast })& \qquad \left\Vert S_{t}^{\ast }f\right\Vert _{q,\kappa
,p}\leq C_{q,\kappa ,p}(S)\left\Vert f\right\Vert _{q+b,\kappa ,p}.
\label{h2'}
\end{align}%
We assume that $C_{q,\kappa ,p}(S)$, $p\in \lbrack 1,\infty ]$, is non
decreasing with respect to $q$ and $\kappa $.
\end{assumption}

We denote 
\begin{align}
C_{q,\kappa ,\infty }(U,S)& =C_{q,\kappa ,\infty }(U)C_{q,\kappa ,\infty
}(S),\quad C_{q,\kappa ,p}(U,S)=C_{q,\kappa ,p}(U)C_{q,\kappa ,p}(S),
\label{h3'} \\
C_{q,\kappa ,\infty ,p}(U,S)& =C_{q,\kappa ,\infty }(U,S)\vee C_{q,\kappa
,p}(U,S).  \label{hh3'}
\end{align}%
Under Assumptions \ref{H1H*1} and \ref{H2H*2}, 
one immediately obtains%
\begin{align}
\left\Vert (S_{t}U_{j})f\right\Vert _{q,-\kappa ,\infty }& \leq C_{q,\kappa
,\infty }(U,S)\left\Vert f\right\Vert _{q+a+b,-\kappa ,\infty },  \label{h}
\\
\left\Vert (S_{t}^{\ast }U_{j}^{\ast })f\right\Vert _{q,\kappa ,p}& \leq
C_{q,\kappa ,p}(U,S)\left\Vert f\right\Vert _{q+a+b,\kappa ,p}.  \label{h'}
\end{align}%
In fact these are the inequalities that we will employ in the following. We
stress that the above constants $C_{q,\kappa ,\infty }(U,S)$ and $%
C_{q,\kappa ,p}(U,S)$ may depend on $a,b$ and are increasing w.r.t. $q$ and $%
\kappa $.

Finally we assume that the (possible) blow-up of $s_{t}\rightarrow \infty $
as $t\rightarrow 0$ is controlled in the following way.

\begin{assumption}
\label{HH3} Let $\theta_0,\lambda>0$ be fixed. We assume that for every $%
\kappa \geq 0$ and $q\in {\mathbb{N}}$ there exist $\pi (q,\kappa )$, $%
\theta _{1}\geq 0$ and $C_{q,\kappa}>0$ such that for every multi-indexes $%
\alpha $ and $\beta $ with $\left\vert \alpha \right\vert +\left\vert \beta
\right\vert \leq q$, $(x,y)\in {\mathbb{R}}^{d}\times {\mathbb{R}}^{d}$ and $%
t\in (0,1]$ one has 
\begin{equation}
(H_{3})\qquad \left\vert \partial _{x}^{\alpha }\partial _{y}^{\beta
}s_{t}(x,y)\right\vert \leq \frac{C_{q,\kappa}}{(\lambda
t)^{\theta_0(q+\theta_1)}}\times \frac{\psi _{\pi (q,\kappa )}(x)}{\psi
_{\kappa }(x-y)}.  \label{h3}
\end{equation}%
We also assume that $\pi (q,k)$ and $C_{q,\kappa}$ are both increasing in $q$
and $\kappa $.
\end{assumption}

This property will be used by means of the following lemma:

\begin{lemma}
\label{lemmaB} Suppose that Assumption \ref{HH3} holds. 

\smallskip

$\mathbf{A.}$ For every $\kappa \geq 0$, $q\in {\mathbb{N}}$ and $p>1$ there
exists $C>0$ such that for every $t\in(0,1]$ and $f$ one has 
\begin{equation}
\left\Vert S_{t}^{\ast }f\right\Vert _{q,\kappa ,p}\leq \frac{C}{(\lambda
t)^{\theta_0(q+\theta_1)}}\left\Vert f\right\Vert _{0,\nu ,1}  \label{B1}
\end{equation}%
where $\nu =\pi (q,\kappa +d)+\kappa +d$

\smallskip

$\mathbf{B.}$ For every $\kappa \geq 0$, $q_1,q_2\in{\mathbb{N}}$, there
exists $C>0$ such that for every $t\in(0,1]$, for every multi-index $\alpha $
with $\left\vert \alpha \right\vert \leq q_{2}$ and $f$ one has 
\begin{equation}
\left\Vert \frac{1}{\psi _{\eta }}S_{t}(\psi _{\kappa }\partial ^{\alpha
}f)\right\Vert _{q_{1},\infty }\leq \frac{C}{(\lambda
t)^{\theta_0(q_1+q_2+\theta_1)}}\,\|f\|_\infty  \label{B2}
\end{equation}%
where $\eta =\pi (q_{1}+q_{2},\kappa +d+1)+\kappa$.
\end{lemma}

\textbf{Proof.} In the sequel, $C$ will denote a positive constant which may
vary from a line to another and which may depend only on $\kappa$ and $q$
for the proof of \textbf{A.} and only on $\kappa, q_1$ and $q_2$ for the
proof of \textbf{B.}

\smallskip

\textbf{A.} Using (\ref{h3}) 
if $\left\vert \alpha \right\vert \leq q,$%
\begin{equation*}
\left\vert \partial ^{\alpha }S_{t}^{\ast }f(x)\right\vert \leq \int
\left\vert \partial _{x}^{\alpha }s_{t}(y,x)\right\vert \times \left\vert
f(y)\right\vert dy\leq \frac{C}{(\lambda t)^{\theta_0(q+\theta_1)}}\int 
\frac{\psi _{\pi (q,\kappa +d)}(y)}{\psi _{\kappa +d}(x-y)}\times \left\vert
f(y)\right\vert dy.
\end{equation*}%
By (\ref{NOT3b}) $\psi _{\kappa +d}(x)/\psi _{\kappa +d}(x-y)\leq C \psi
_{\kappa +d}(y)$ so that%
\begin{align*}
\psi _{\kappa +d}(x)\left\vert \partial ^{\alpha }S_{t}^{\ast
}f(x)\right\vert &\leq \frac{C}{(\lambda t)^{\theta_0(q+\theta_1)}} \int 
\frac{ \psi _{\kappa +d}(x)\psi _{\pi (q,\kappa +d)}(y)}{\psi _{\kappa
+d}(x-y)} \times \left\vert f(y)\right\vert dy \\
&\leq \frac{C}{(\lambda t)^{\theta_0(q+\theta_1)}} \int \psi _{\pi (q,\kappa
+d)+\kappa +d}(y)\times \left\vert f(y)\right\vert dy \\
&=\frac{C}{(\lambda t)^{\theta_0(q+\theta_1)}}\left\Vert f\right\Vert
_{0,\nu ,1}.
\end{align*}%
We conclude that 
\begin{equation*}
\left\Vert S_{t}^{\ast }f\right\Vert _{q,\kappa +d,\infty }\leq \frac{C}{%
(\lambda t)^{\theta_0(q+\theta_1)}}\left\Vert f\right\Vert _{0,\nu ,1}.
\end{equation*}
By (\ref{NOT5a}) $\left\Vert S_{t}^{\ast }f\right\Vert _{q,\kappa ,p}\leq
C\left\Vert S_{t}^{\ast }f\right\Vert _{q,\kappa +d,\infty }$ so the proof
of (\ref{B1}) is completed.

\smallskip

\textbf{B.} Let $\gamma $ with $\left\vert \gamma \right\vert \leq q_{1}$.
Using integration by parts%
\begin{eqnarray*}
\partial ^{\gamma }S_{t}(\psi _{\kappa }\partial ^{\alpha }f)(x) &=&\int_{{%
\mathbb{R}}^{d}}\partial _{x}^{\gamma }s_{t}(x,y)\psi _{\kappa }(y)\partial
^{\alpha }f(y)dy \\
&=&(-1)^{\left\vert \alpha \right\vert }\int_{{\mathbb{R}}^{d}}\partial
_{y}^{\alpha }(\partial _{x}^{\gamma }s_{t}(x,y)\psi _{\kappa }(y))\times
f(y)dy.
\end{eqnarray*}%
Using (\ref{NOT3c}), (\ref{h3}) 
and (\ref{NOT3b}), 
it follows that 
\begin{align*}
\left\vert \partial ^{\gamma }S_{t}(\psi _{\kappa }\partial ^{\alpha
}f)(x)\right\vert &\leq \int_{{\mathbb{R}}^{d}}\left\vert \partial
_{y}^{\alpha }(\partial _{x}^{\gamma }s_{t}(x,y)\psi _{\kappa
}(y))\right\vert \times \left\vert f(y)\right\vert dy \\
&\leq \int_{{\mathbb{R}}^{d}}\left\vert s_{t}(x,y)\psi _{\kappa
}(y)\right\vert _{q_{1}+q_{2}}\times \left\vert f(y)\right\vert dy \\
&\leq C \int_{{\mathbb{R}}^{d}}\left\vert s_{t}(x,y)\right\vert
_{q_{1}+q_{2}}\psi _{\kappa }(y)\times \left\vert f(y)\right\vert dy \\
&\leq \frac{C}{(\lambda t)^{\theta_0(q_1+q_2+\theta_1)}} \left\Vert
f\right\Vert _{\infty }\int_{{\mathbb{R}}^{d}}\frac{\psi _{\pi
(q_{1}+q_{2},\kappa +d+1)}(x)}{\psi _{\kappa +d+1}(x-y)}\times \psi _{\kappa
}(y)dy \\
&\leq \frac{C}{(\lambda t)^{\theta_0(q_1+q_2+\theta_1)}} \left\Vert
f\right\Vert _{\infty }\int_{{\mathbb{R}}^{d}}\frac{\psi _{\pi
(q_{1}+q_{2},\kappa +d+1)+\kappa }(x)}{\psi _{d+1}(x-y)}dy \\
&\leq \frac{C}{(\lambda t)^{\theta_0(q_1+q_2+\theta_1)}}\left\Vert
f\right\Vert _{\infty }\psi _{\pi (q_{1}+q_{2},\kappa +d+1)+\kappa }(x).
\end{align*}%
%
%
%
%
%
%
%
%
%
%
%
%
This implies (\ref{B2}). $\square $

\medskip

We are now able to give the \textquotedblleft regularity
lemma\textquotedblright . This is the core of our approach.

\begin{lemma}
\label{Reg} Suppose that Assumption \ref{H1H*1}, \ref{H2H*2} and \ref{HH3}
hold. 
We fix $t\in (0,1]$, $m\geq 1$ and $\delta_{i}>0$, $i=1,\ldots ,m$ such that 
$\sum_{i=1}^{m}\delta_{i}=t.$

\smallskip

$\mathbf{A.}$ There exists a function $\tilde{p}_{\delta_{1},...,\delta_{m}}%
\in C^{\infty }({\mathbb{R}}^{d}\times {\mathbb{R}}^{d})$ such that 
\begin{equation}
\prod_{i=1}^{m-1}(S_{\delta_{i}}U_{i})S_{\delta_{m}}f(x)=\int \tilde{p}%
_{\delta_{1},...,\delta_{m}}(x,y)f(y)dy.  \label{h6}
\end{equation}

$\mathbf{B.}$ We fix $q_{1},q_{2}\in {\mathbb{N}},\kappa \geq 0,p>1$ and we
denote $q=q_{1}+q_{2}+(a+b)(m-1).$ One may find universal constants $C,\chi ,%
\bar{p}\geq 1$ (depending on $\kappa ,p$ and $q_{1}+q_{2})$ such that for
every multi-index $\beta $ with $\left\vert \beta \right\vert \leq q_{2}$
and every $x\in {\mathbb{R}}^{d}$%
\begin{equation}
\left\Vert \partial _{x}^{\beta }\tilde{p}_{\delta _{1},...,\delta
_{m}}(x,\cdot )\right\Vert _{q_{1},\kappa ,p}\leq C\Big(\frac{2m}{\lambda t}%
\Big)^{\theta _{0}(q_{1}+q_{2}+d+2\theta _{1})}\Big(C_{q,\chi ,\bar{p}%
,\infty }(U,S)\Big(\frac{2m}{\lambda t}\Big)^{\theta _{0}(a+b)}\Big)%
^{m-1}\psi _{\chi }(x).  \label{h7}
\end{equation}%
%
%
%
%
%
%
%
%
%
%
%
%
\end{lemma}

\textbf{Proof.} \textbf{A.} For $g=g(x,y)$, we denote $g^{x}(y):=g(x,y)$. By
the very definition of $U_{i}^{\ast }$ one has%
\begin{equation*}
S_{t}U_{i}f(x)=\int_{{\mathbb{R}}^{d}}U_{i}^{\ast }s_{t}^{x}(y)f(y)dy.
\end{equation*}%
As a consequence, one gets the kernel in (\ref{h6}): 
\begin{equation*}
\tilde{p}_{\delta_{1},\ldots ,\delta_{m}}(x,y)=\int_{{\mathbb{R}}^{d\times
(m-1)}}U_{1}^{\ast }s_{\delta_{1}}^{x}(y_{1})\Big(\prod_{j=2}^{m-1}U_{j}^{%
\ast }s_{\delta_{j}}^{y_{j-1}}(y_{j})\Big)s_{\delta_{m}}(y_{m-1},y)dy_{1}%
\cdots dy_{m-1},
\end{equation*}%
and the regularity immediately follows.

\smallskip

\textbf{B.} We split the proof in several steps.

\smallskip \textbf{Step 1: decomposition}. Since $\sum_{i=1}^{m}\delta_{i}=t$
we may find $j\in \{1,...,m\}$ such that $\delta_{j}\geq \frac{t}{m}.$ We
fix this $j$ and we write%
\begin{equation*}
\prod_{i=1}^{m-1}(S_{\delta_{i}}U_{i})S_{\delta_{m}}=Q_{1}Q_{2}
\end{equation*}%
with%
\begin{equation*}
Q_{1}=\prod_{i=1}^{j-1}(S_{\delta_{i}}U_{i})S_{\frac{1}{2}\delta_{j}}\quad %
\mbox{and}\quad Q_{2}=S_{\frac{1}{2}\delta_{j}}U_{j}\prod_{i=j+1}^{m-1}(S_{%
\delta_{i}}U_{i})S_{\delta_{m}}=S_{\frac{1}{2}\delta_{j}}%
\prod_{i=j}^{m-1}(U_{i}S_{\delta_{i+1}}).
\end{equation*}%
Here we use the semi-group property $S_{\frac{1}{2}\delta_{j}}S_{\frac{1}{2}%
\delta_{j}}=S_{\delta_{j}}.$

We suppose that $j\leq m-1$. In the case $j=m$ the proof is analogous but
simpler. We will use Lemma \ref{lemmaB} in order to estimate the terms
corresponding to each of these two operators. As already seen, both $Q_1$
and $Q_2$ are given by means of smooth kernels, that we call $p_1(x,y)$ and $%
p_2(x,y)$ respectively.

\smallskip

\textbf{Step 2}. We take $\beta $ with $\left\vert \beta \right\vert \leq
q_{2}$ and we denote $g^{\beta ,x}(y):=\partial _{x}^{\beta }g(x,y)$. For $%
h\in L^{1}$ we write%
\begin{align*}
& \int_{{\mathbb{R}}^{d}}h(z)\partial _{x}^{\beta }\tilde{p}_{\delta
_{1},...,\delta _{m}}(x,z)dz=\int_{{\mathbb{R}}^{d}}h(z)\int_{{\mathbb{R}}%
^{d}}\partial _{x}^{\beta }p_{1}(x,y)p_{2}(y,z)dydz \\
& \qquad =\int_{{\mathbb{R}}^{d}}\partial _{x}^{\beta }p_{1}(x,y)\int
h(z)p_{2}(y,z)dzdy=\int_{{\mathbb{R}}^{d}}\partial _{x}^{\beta
}p_{1}(x,y)Q_{2}h(y)dy \\
& \qquad =\int_{{\mathbb{R}}^{d}}Q_{2}^{\ast }p_{1}^{\beta ,x}(y)h(y)dy.
\end{align*}%
It follows that 
\begin{equation*}
\partial _{x}^{\beta }\tilde{p}_{\delta _{1},...,\delta
_{m}}(x,z)=Q_{2}^{\ast }p_{1}^{\beta ,x}(z)=\prod_{i=1}^{m-j}(S_{\delta
_{m-i+1}}^{\ast }U_{m-i}^{\ast })S_{\frac{1}{2}\delta _{j}}^{\ast
}p_{1}^{\beta ,x}(z).
\end{equation*}%
We will use (\ref{h'}) $m-j$ times first and (\ref{B1}) then. We denote 
\begin{equation*}
q_{1}^{\prime }=q_{1}+(m-j)(a+b)
\end{equation*}%
and we write 
\begin{equation}
\begin{array}{rl}
\Vert \partial _{x}^{\beta }\tilde{p}_{\delta _{1},...,\delta _{m}}(x,\cdot
)\Vert _{q_{1},\kappa ,p} & \leq C_{q_{1}^{\prime },\kappa
,p}^{m-j}(U,S)\Vert S_{\frac{1}{2}\delta _{j}}^{\ast }p_{1}^{\beta ,x}\Vert
_{q_{1}^{\prime },\kappa ,p}\smallskip \\ 
& \displaystyle\leq C_{q_{1}^{\prime },\kappa ,p}^{m-j}(U,S)\,C\Big(\frac{2m%
}{\lambda t}\Big)^{\theta _{0}(q_{1}^{\prime }+\theta _{1})}\Vert
p_{1}^{\beta ,x}\Vert _{0,\nu ,1}%
\end{array}
\label{h9}
\end{equation}%
with%
\begin{equation*}
\nu =\pi (q_{1}^{\prime },\kappa +d)+\kappa +d.
\end{equation*}%
\textbf{Step 3.} We denote $g_{z}(u)=\prod_{l=1}^{d}1_{(0,\infty
)}(u_{l}-z_{l})$, so that $\delta _{0}(u-z)=\partial _{u}^{\rho }g_{z}(u)$
with $\rho =(1,2,\ldots ,d).$ We take $\mu =\nu +d+1$ and we formally write 
\begin{equation*}
p_{1}(x,z)=\frac{1}{\psi _{\mu }(z)}Q_{1}(\psi _{\mu }\partial ^{\rho
}g_{z})(x).
\end{equation*}%
%
%
%
%
%
%
%
%
%
%
This formal equality can be rigorously written by using the regularization
by convolution of the Dirac function.

We denote%
\begin{equation*}
q_{2}^{\prime }=q_{2}+(j-1)(a+b),\quad \eta =\pi (d+q_{2}^{\prime },\mu
+d+1)+\mu
\end{equation*}%
and we write%
\begin{equation*}
|p_{1}^{\beta ,x}(z)|=|\partial _{x}^{\beta }p_{1}(x,z)|\leq \frac{\psi
_{\eta }(x)}{\psi _{\mu }(z)}\Big\Vert\frac{1}{\psi _{\eta }}\partial
^{\beta }Q_{1}(\psi _{\mu }\partial ^{\rho }g_{z})\Big\Vert_{\infty }.
\end{equation*}%
Since $\mu =\nu +d+1$, $\int \psi _{\nu }\times \frac{1}{\psi _{\mu }}%
<\infty $, so using (\ref{NOT3c}), we obtain (recall that $\left\vert \beta
\right\vert \leq q_{2})$ 
\begin{align*}
\Vert p_{1}^{\beta ,x}\Vert _{0,\nu ,1}& \leq C\psi _{\eta }(x)\sup_{z\in {%
\mathbb{R}}^{d}}\big\Vert\frac{1}{\psi _{\eta }}\partial ^{\beta }Q_{1}(\psi
_{\mu }\partial ^{\rho }g_{z})\big\Vert_{\infty }\leq C\psi _{\eta
}(x)\sup_{z\in {\mathbb{R}}^{d}}\big\Vert\frac{1}{\psi _{\eta }}Q_{1}(\psi
_{\mu }\partial ^{\rho }g_{z})\big\Vert_{q_{2},\infty } \\
& \leq C\psi _{\eta ^{\prime }}(x)\sup_{z\in {\mathbb{R}}^{d}}\big\Vert %
Q_{1}(\psi _{\mu }\partial ^{\rho }g_{z})\big\Vert_{q_{2},-\eta ,\infty }.
\end{align*}%
Using (\ref{h}) $j-1$ times and (\ref{B2}) (with $\kappa =\mu )$ we get 
\begin{eqnarray*}
\left\Vert Q_{1}(\psi _{\mu }\partial ^{\rho }g_{z})\right\Vert
_{q_{2},-\eta ,\infty } &\leq &C_{q_{2}^{\prime },\eta ,\infty
}^{j-1}(U,S)\Vert S_{\frac{1}{2}\delta _{j}}(\psi _{\mu }\partial ^{\rho
}g_{z})\Vert _{q_{2}^{\prime },-\eta ,\infty } \\
&\leq &C_{q_{2}^{\prime },\eta ,\infty }^{j-1}(U,S)\left\Vert
g_{z}\right\Vert _{\infty }\,C\Big(\frac{2m}{\lambda t}\Big)^{\theta
_{0}(q_{2}^{\prime }+d+\theta _{1})}.
\end{eqnarray*}%
Since $\left\Vert g_{z}\right\Vert _{\infty }=1$ we obtain%
\begin{equation*}
\Vert p_{1}^{\beta ,x}\Vert _{0,\nu ,1}\leq \psi _{\eta }(x)C_{q_{2}^{\prime
},\eta ,\infty }^{j-1}(U,S)\,C\Big(\frac{2m}{\lambda t}\Big)^{\theta
_{0}(q_{2}^{\prime }+d+\theta _{1})}.
\end{equation*}%
By inserting in (\ref{h9}) we obtain (\ref{h7}), so the proof is completed. $%
\square $

\section{Proofs of the main results}

\label{sect:proofs}

In this section we prove Theorem \ref{TransferBIS}. But before we give an
intermediary result, Theorem \ref{Transfer} below, which is more precise
concerning constants. Let us introduce some notation. For $\delta \geq 0$ we
denote%
\begin{equation}
\Phi _{n}(\delta )=\varepsilon _{n}\Lambda _{n}\times \lambda _{n}^{-\theta
_{0}(a+b+\delta )}.  \label{R7'}
\end{equation}%
We recall that the constants $\varepsilon_n$, $a$, $\Lambda_n$, $b$ and $%
\lambda_n$ are defined in Assumption \ref{A1A*1}, Assumption \ref{A2A*2} and
Assumption \ref{A3}. 
Under Assumption \ref{A3}, $\lambda_n\leq\gamma\lambda_{n+1}$ for some $%
\gamma\geq 1$, so we have 
\begin{equation}
\Phi _{n}(\delta )\leq \gamma ^{1+\theta _{0}(a+b+\delta )}\Phi
_{n+1}(\delta ).  \label{TRa}
\end{equation}
For $\kappa \geq 0,\eta \geq 0$ we set 
\begin{equation}
\Psi _{\eta ,\kappa }(x,y):=\frac{\psi _{\kappa }(y)}{\psi _{\eta }(x)}%
,\quad (x,y)\in {\mathbb{R}}^{d}\times {\mathbb{R}}^{d}.  \label{R7''}
\end{equation}

Our intermediary result concerning the regularity of the semigroup $%
(P_t)_{t\geq 0}$ is the following.

\begin{theorem}
\label{Transfer} Suppose that Assumption \ref{A1A*1}, \ref{A2A*2} , \ref{A3}
and \ref{A5} 
hold. Moreover we suppose there exists $\delta >0$ such that 
\begin{equation}
\limsup_{n\to\infty}\Phi _{n}(\delta )<\infty ,  \label{TR6}
\end{equation}%
$\Phi_n(\delta)$ being given in (\ref{R7'}). Then the following statements
hold.

\smallskip

\textbf{A}. $P_{t}f(x)=\int_{{\mathbb{R}}^{d}}p_{t}(x,y)dy$ with $p_{t}\in
C^{\infty }({\mathbb{R}}^{d}\times {\mathbb{R}}^{d}).$

\smallskip

\textbf{B}. Let $n_\ast\in {\mathbb{N}}$ and $\delta _{\ast }>0$ be such that \ 
\begin{equation}
\overline{\Phi }_{\ast}:=\sup_{n\geq n_\ast}\Phi
_{n}(\delta _{\ast })<\infty .  \label{TR6a}
\end{equation}%
We fix $q\in {\mathbb{N}}$, $p>1$, $\varepsilon _{\ast }>0,$ $\kappa \geq 0$
and we put $\mathfrak{m}=1+\frac{q+2d/p_{\ast }}{\delta _{\ast }}$ with $%
p_{\ast }$ the conjugate of $p$. 
There exist $C\geq 1$ and $\eta >1$ (depending on $q,p,\varepsilon _{\ast
},\delta _{\ast },\kappa $ and $\gamma $) such that for every $t\in(0,1]$ 
\begin{eqnarray}
\left\Vert \Psi _{\eta ,\kappa }p_{t}\right\Vert _{q,p} &\leq &C\times
Q_\ast(q,\mathfrak{m})\times t^{-\theta _{0}((a+b)\mathfrak{m}+q+2d/p_{\ast
})(1+\varepsilon _{\ast })}\quad with  \label{TR6'} \\
Q_\ast(q,\mathfrak{m})&=&\Big(\frac{1}{\lambda _{n_\ast}^{\theta _{0}(a+b)%
\mathfrak{m}+q+2d/p_{\ast }}}+\overline{\Phi }_\ast^{\mathfrak{m}}\Big)^{1+\varepsilon _{\ast }}.  \label{TR6''}
\end{eqnarray}

\smallskip

\textbf{C}. Let $p>2d.$ Set $\bar{\mathfrak{m}}=1+\frac{q+1+2d/p_{\ast }}{%
\delta _{\ast }}$. There exist $C\geq 1,\eta > 1$ (depending on $%
q,p,\varepsilon _{\ast },\delta _{\ast },\kappa $) such that for every $t\in (0,1]$, 
 $x,y\in {\mathbb{R}}^{d}$ and for every multi-indexes $\alpha ,\beta $ such that 
$\left\vert \alpha \right\vert +\left\vert \beta \right\vert \leq q$, 
\begin{equation}
\left\vert \partial _{x}^{\alpha }\partial _{y}^{\beta
}p_{t}(x,y)\right\vert \leq C\times Q_\ast(q+1,\bar{\mathfrak{m}})\times
t^{-\theta _{0}((a+b)\bar{\mathfrak{m}}+q+1+2d/p_{\ast })(1+\varepsilon
_{\ast })}\times \frac{\psi _{\eta +\kappa }(x)}{\psi _{\kappa }(x-y)}
\label{TR6d}
\end{equation}%
%
%
%
%
%
%
%
%
%
%
%
\end{theorem}

\begin{remark}
We stress that in hypothesis (\ref{TR6a}) the order of derivation $q$ does
not appear. However the conclusions (\ref{TR6'}) and (\ref{TR6d}) hold for
every $q.$ The motivation of this\ is given by the following heuristics. The
hypothesis (\ref{TR5}) says that the semi-group $P_{t}^{n}$ has a
regularization effect controlled by $1/(\lambda _{n}t)^{\theta _{0}}.$ If we
want to decouple this effect $m_{0}$ times we write $%
P_{t}^{n}=P_{t/m_{0}}^{n}....P_{t/m_{0}}^{n}$ and then each of the $m_{0}$
operators $P_{t/m_{0}}^{n}$ acts with a regularization effect of order $%
(\lambda _{n}\times t/m_{0})^{\theta _{0}}.$ But this heuristics does not
work directly: in order to use it, in the proof we have to develop a Taylor
expansion coupled with the interpolation criterion studied in Section \ref%
{sect:reg}.
\end{remark}

\textbf{Proof. Step 0: constants and parameters set-up.} We first choose
some parameters which will be used in the following steps. To begin we
stress that we work with measures on ${\mathbb{R}}^{d}\times {\mathbb{R}}^{d}
$ so the dimension of the space is $2d$ (and not $d).$ We recall that in our
statement the quantities $q,d,p,\delta _{\ast },\varepsilon _{\ast },\kappa $
and $n$ are given and fixed. In the following we will denote by $C $ a
constant depending on all these parameters and which may change from a line
to another. We define 
\begin{equation}
m_{0}=1+\Big\lfloor\frac{q+2d/p_{\ast }}{\delta _{\ast }}\Big\rfloor>0
\label{H4}
\end{equation}%
and given $h\in {\mathbb{N}}$ we denote 
\begin{equation}
\rho _{h}=\frac{(a+b)m_{0}+q+2d/p_{\ast }}{2h}.  \label{H5'}
\end{equation}%
Notice that this is equal to the constant $\rho _{h}$ defined in (\ref{reg5}%
) corresponding to $k=(a+b)m_{0}$ and $q$ and to $2d$ (instead of $d).$

\smallskip

\textbf{Step 1: a Lindeberg-type method to decompose $P_t-P^n_t$.} We fix
(once for all) $t\in(0,1]$ and we write%
\begin{equation*}
P_{t}f-P_{t}^{n}f=\int_{0}^{t}\partial
_{s}(P_{t-s}^{n}P_{s})fds=\int_{0}^{t}P^n
_{t-s}(L-L_{n})P_{s}fds=\int_{0}^{t}P^n_{t-s}\Delta _{n}P_{s}fds.
\end{equation*}%
We iterate this formula $m_{0}$ times (with $m_{0}$ chosen in (\ref{H4}))
and we obtain%
\begin{equation}
P_{t}f(x)-P_{t}^{n}f(x)=\sum_{m=1}^{m_{0}-1}I_{n}^{m}f(x)+R_{n}^{m_{0}}f(x)
\label{R2}
\end{equation}%
with (we put $t_{0}=t)$%
\begin{align*}
I_{n}^{m}f(x)&=\int_{0}^{t}dt_{1}\int_{0}^{t_{1}}dt_{2}...%
\int_{0}^{t_{m-1}}dt_{m}\prod_{i=0}^{m-1}(P_{t_{i}-t_{i+1}}^{n}\Delta
_{n})P_{t_{m}}^{n}f(x),\quad 1\leq m\leq m_0-1, \\
R_{n}^{m_{0}}f(x)&=\int_{0}^{t}dt_{1}\int_{0}^{t_{1}}dt_{2}...%
\int_{0}^{t_{m_{0}-1}}dt_{m_{0}}\prod_{i=0}^{m_{0}-1}(P_{t_{i}-t_{i+1}}^{n}%
\Delta _{n})P_{t_{m_{0}}}f(x).
\end{align*}

In order to analyze $I_{n}^{m}f$ we use Lemma \ref{Reg} for the semigroup $%
S_{t}=P_{t}^{n}$ and for the operators $U_{i}=\Delta _{n}=L-L_{n}$ (the same
for each $i$), with $\delta _{i}=t_{i}-t_{i+1}$, $i=0,\ldots ,m$ (with $%
t_{m+1}=0$). So the hypotheses 
(\ref{h1}) and (\ref{h1'}) in Assumption \ref{H1H*1} coincide with the
requests 
(\ref{TR3}) and (\ref{TR3'}) in Assumption \ref{A1A*1}. And we have $%
C_{q,\kappa ,\infty }(U)=C_{q,\kappa ,p}(U)=C\varepsilon _{n}.$ Moreover the
hypotheses 
(\ref{h2}) and (\ref{h2'}) in Assumption \ref{H2H*2} coincide with the
hypotheses 
(\ref{TR2}) and (\ref{TR2'}) in Assumption \ref{A2A*2}. And we have $%
C_{q,\kappa ,\infty }(P^{n})=C_{q,\kappa ,p}(P^{n})=\Lambda _{n}$. Hence, 
\begin{equation}
C_{q,\kappa ,\infty ,p}(\Delta _{n},P^{n})=C\,\varepsilon _{n}\times \Lambda
_{n}.  \label{app1}
\end{equation}%
Finally, the hypothesis 
(\ref{h3}) in Assumption \ref{HH3} coincides with 
(\ref{TR5}) in Assumption \ref{A3}. So, we can apply Lemma \ref{Reg}: by
using (\ref{h6}) we obtain%
\begin{equation*}
I_{n}^{m}f(x)=\int_{0}^{t}dt_{1}...\int_{0}^{t_{m-1}}dt_{m}\int
p_{t-t_{1},t_{1}-t_{2},...,t_{m}}^{n,m}(x,y)f(y)dy.
\end{equation*}%
We denote%
\begin{equation*}
\phi
_{t}^{n,m_{0}}(x,y)=p_{t}^{n}(x,y)+\sum_{m=1}^{m_{0}-1}\int_{0}^{t}dt_{1}...%
\int_{0}^{t_{m-1}}dt_{m}p_{t-t_{1},t_{2}-t_{1},...,t_{m}}^{n,m}(x,y)
\end{equation*}%
so that (\ref{R2}) reads%
\begin{equation*}
\int f(y)P_{t}(x,dy)=\int f(y)\phi _{t}^{n,m_{0}}(x,y)dy+R_{n}^{m_{0}}f(x).
\end{equation*}%
We recall that $\Psi _{\eta ,\kappa }$ is defined in (\ref{R7''}) and we
define the measures on ${\mathbb{R}}^{d}\times {\mathbb{R}}^{d}$ defined by 
\begin{equation*}
\mu ^{\eta ,\kappa }(dx,dy)=\Psi _{\eta ,\kappa }(x,y)P_{t}(x,dy)dx\quad %
\mbox{and}\quad \mu _{n}^{\eta ,\kappa ,m_{0}}(dx,dy)=\Psi _{\eta ,\kappa
}(x,y)\phi _{t}^{n,m_{0}}(x,y)dxdy.
\end{equation*}%
So, the proof consists in applying Lemma \ref{REG} to $\mu =\mu ^{\eta
,\kappa }$ and $\mu _{n}=\mu _{n}^{\eta ,\kappa ,m_{0}}$.

\smallskip

\textbf{Step 2: analysis of the principal term.} We study here the estimates
for $f_{n}(x,y)=\Psi _{\eta ,\kappa }\phi _{t}^{n,m_{0}}(x,y)$ which are
required in (\ref{reg9}).

We first use (\ref{h7}) in order to get estimates for $%
p_{t-t_{1},t_{2}-t_{1},...,t_{m}}^{n,m}(x,y)$. We fix $q_{1},q_{2}\in {\ 
\mathbb{N}},\kappa \geq 0,p>1$ and we recall that in Lemma \ref{Reg} we
introduced $\overline{q}=q_{1}+q_{2}+(a+b)(m_{0}-1).$ Moreover in Lemma \ref%
{Reg} one produces $\chi $ such that (\ref{h7}) holds true: for every
multi-index $\beta $ with $\left\vert \beta \right\vert \leq q_{2}$%
\begin{equation*}
\begin{array}{l}
\left\Vert \psi _{\kappa }\partial _{x}^{\beta
}p_{t-t_{1},t_{1}-t_{2},...,t_{m}}^{n,m}(x,\cdot )\right\Vert
_{q_{1},p}\smallskip \\ 
\displaystyle\quad \leq C\Big(\frac{1}{\lambda _{n}t}\Big)^{\theta
_{0}(q_{1}+q_{2}+d+2\theta _{1})}\times \left( \varepsilon _{n}\Lambda _{n}%
\Big(\frac{1}{\lambda _{n}t}\Big)^{\theta _{0}(a+b)}\right) ^{m}\psi _{\chi
}(x).%
\end{array}%
\end{equation*}%
%
Denote%
\begin{equation*}
\xi _{1}(q)=q+d+2\theta _{1}+m_{0}(a+b),\qquad \omega _{1}(q)=q+d+2\theta
_{1}.
\end{equation*}%
With this notation, if $\left\vert \beta \right\vert \leq q_{2}$ we have 
\begin{eqnarray}
\left\Vert \psi _{\kappa }\partial _{x}^{\beta }\phi _{t}^{n,m_{0}}(x,\cdot
)\right\Vert _{q_{1},p} &\leq &C\Big(\frac{1}{\lambda _{n}t}\Big)^{\theta
_{0}(q_{1}+q_{2}+d+2\theta _{1})}\times \left( \varepsilon _{n}\Lambda _{n}%
\Big(\frac{1}{\lambda _{n}t}\Big)^{\theta _{0}(a+b)}\right) ^{m_{0}}\psi
_{\chi }(x)  \label{h12} \\
&=&Ct^{-\theta _{0}\xi _{1}(q_{1}+q_{2})}\lambda _{n}^{-\theta _{0}\omega
_{1}(q_{1}+q_{2})}\Phi _{n}^{m_{0}}(0)\psi _{\chi }(x),
\end{eqnarray}
where $\Phi _{n}(\delta )$ is the constant defined in (\ref{R7'}). We take $%
l=2h+q,l^{\prime }=2h$ and we take $q(l)=l+(a+b)m_{0}.$ Moreover we fix $%
q_{1}$ and $q_{2}$ (so $q=q_{1}+q_{2}\leq l)$ and we take $\chi $ to be the
one in (\ref{h12}). Moreover we take $\eta $ sufficiently large in order to
have $p\eta -2h-p\chi \geq d+1.$ This guarantees that 
\begin{equation}
\int_{{\mathbb{R}}^{d}}\frac{dx}{\psi _{p\eta -l^{\prime }-p\chi }(x)}%
=C<\infty .  \label{h14}
\end{equation}%
By (\ref{NOT3c}) and (\ref{h12}) 
\begin{align*}
& \left\Vert \Psi _{\eta ,\kappa }\phi _{t}^{n,m_{0}}\right\Vert
_{l,l^{\prime },p}^{p}\leq C\sum_{\left\vert \alpha \right\vert +\left\vert
\beta \right\vert \leq l}\int_{{\mathbb{R}}^{d}}\int_{{\mathbb{R}}^{d}}\Psi
_{\eta ,\kappa }^{p}(x,y)\left\vert \partial _{x}^{\alpha }\partial
_{y}^{\beta }\phi _{t}^{n,m_{0}}(x,y)\right\vert ^{p}\psi _{l^{\prime
}}(x)\psi _{l^{\prime }}(y)dydx \\
& \quad =C\sum_{\left\vert \alpha \right\vert +\left\vert \beta \right\vert
\leq l}\int_{{\mathbb{R}}^{d}}\frac{1}{\psi _{p\eta -l^{\prime }}(x)}\int_{{%
\ \mathbb{R}}^{d}}\left\vert \psi _{\kappa +l^{\prime }/p}(y)\partial
_{x}^{\alpha }\partial _{y}^{\beta }\phi _{t}^{n,m_{0}}(x,y)\right\vert
^{p}dydx \\
& \quad \leq C\sum_{\left\vert \alpha \right\vert +\left\vert \beta
\right\vert \leq l}\int_{{\mathbb{R}}^{d}}\frac{1}{\psi _{p\eta -l^{\prime
}}(x)}\left\Vert \psi _{\kappa +l^{\prime }/p}\partial _{x}^{\alpha }\phi
_{t}^{n,m_{0}}(x,\cdot )\right\Vert _{\left\vert \beta \right\vert ,p}^{p}dx
\\
& \quad \leq C(t^{-\theta _{0}\xi _{1}(l)}\lambda _{n}^{-\theta _{0}\omega
_{1}(l)}\Phi _{n}(0))^{pm_{0}}\int_{{\mathbb{R}}^{d}}\frac{dx}{\psi _{p\eta
-l^{\prime }-p\chi }(x)}.
\end{align*}%
We conclude that 
\begin{equation}
\left\Vert \Psi _{\eta ,\kappa }\phi _{t}^{n,m_{0}}\right\Vert
_{2h+q,2h,p}\leq Ct^{-\theta _{0}\xi _{1}(q+2h)}\times \lambda _{n}^{-\theta
_{0}\omega _{1}(q+2h)}\Phi _{n}^{m_{0}}(0)=:\theta (n).  \label{R6}
\end{equation}%
By (\ref{TRa}) $\theta (n)\uparrow +\infty $ and $\Theta \theta (n)\geq
\theta (n+1)$ with 
\begin{equation*}
\Theta =\gamma ^{\theta _{0}((a+b)m_{0}+q+2h+d+2\theta _{1})+m_{0}}\geq 1.
\end{equation*}%
In the following we will choose $h$ sufficiently large, depending on $\delta
_{\ast },m_{0},q,d$ and $p.$ So $\Theta $ is a constant depending on $\delta
_{\ast },m_{0},q,d,a,b$,$\gamma $ and $p,$ as the constants considered in
the statement of our theorem.

\smallskip

\textbf{Step 3: analysis of the remainder}. We study here $%
d_{m_{0}}(n):=d_{(a+b)m_{0}}(\mu ^{\eta ,\kappa },\mu _{n}^{\eta ,\kappa
,m_{0}})$ as required in (\ref{reg10}): we prove that, if $\eta \geq \kappa
+d+1,$ then%
\begin{equation}
d_{m_{0}}(n)\leq C(\Lambda _{n}\varepsilon _{n})^{m_{0}}\leq \lambda
_{n}^{\theta _{0}(a+b+\delta _{\ast })m_{0}}\Phi _{n}^{m_{0}}(\delta _{\ast
}).  \label{R9}
\end{equation}

Using first $(A_{1})$ and $(A_{2})$ (see (\ref{TR3}) and (\ref{TR2})) and
then $(A_{4})$ (see (\ref{R7})) we obtain 
\begin{equation*}
\left\Vert \prod_{i=0}^{m_{0}-1}(P_{t_{i}-t_{i+1}}^{n}\Delta
_{n})P_{t_{m_{0}}}f\right\Vert _{0,-\kappa ,\infty }\leq C\left\Vert
f\right\Vert _{(a+b)m_{0},-\kappa ,\infty }(\Lambda _{n}\varepsilon
_{n})^{m_{0}}
\end{equation*}%
which gives 
\begin{equation*}
\left\Vert R_{n}^{m_{0}}f\right\Vert _{0,-\kappa ,\infty }\leq C\left\Vert
f\right\Vert _{(a+b)m_{0},-\kappa ,\infty }(\Lambda _{n}\varepsilon
_{n})^{m_{0}}.
\end{equation*}%
Using now the equivalence between (\ref{NOT6a}) and (\ref{NOT6b}) we obtain%
\begin{equation}
\left\Vert \frac{1}{\psi _{\kappa }}R_{n}^{m_{0}}(\psi _{\kappa
}f)\right\Vert _{\infty }\leq C\left\Vert f\right\Vert _{(a+b)m_{0},\infty
}(\Lambda _{n}\varepsilon _{n})^{m_{0}}.  \label{R8}
\end{equation}%
We take now $g\in C^{\infty }({\mathbb{R}}^{d}\times {\mathbb{R}}^{d}),$ we
denote $g_{x}(y)=g(x,y),$ and we write 
\begin{align*}
& \left\vert \int_{{\mathbb{R}}^{d}\times {\mathbb{R}}^{d}}g(x,y)(\mu ^{\eta
,\kappa }-\mu _{n}^{,\kappa ,m_{0}})(dx,dy)\right\vert \\
& \quad \leq \int_{{\mathbb{\ R}}^{d}}\frac{dx}{\psi _{\eta }(x)}\left\vert
\int_{{\mathbb{R}}^{d}}g_{x}(y)\psi _{\kappa }(y)(P_{t}(x,dy)-\phi
_{t}^{n,m_{0}}(x,y))dy\right\vert \\
& \quad \leq \int_{{\mathbb{R}}^{d}}\frac{dx}{\psi _{\eta -\kappa }(x)}%
\left\vert \frac{1}{\psi _{\kappa }(x)}R_{n}^{m_{0}}(\psi _{\kappa
}g_{x})(x)\right\vert dx \\
& \quad \leq C\sup_{x\in {\mathbb{R}}^{d}}\left\Vert g_{x}\right\Vert
_{(a+b)m_{0},\infty }(\Lambda _{n}\varepsilon _{n})^{m_{0}}
\end{align*}%
the last inequality being a consequence of (\ref{R8}) and of $\eta -\kappa
\geq d+1$. Now (\ref{R9}) is proved because $\sup_{x\in {\mathbb{R}}%
^{d}}\left\Vert g_{x}\right\Vert _{(a+b)m_{0},\infty }$ $\leq \left\Vert
g\right\Vert _{(a+b)m_{0},\infty }$.

\smallskip

\textbf{Step 4: use of Lemma \ref{REG} and proof of A. and B.} We recall
that $\rho _{h}$ is defined in (\ref{H5'}) and we estimate 
\begin{equation*}
d_{m_{0}}(n)\times \theta (n)^{\rho _{h}}\leq Ct^{-\theta _{0}\xi
_{2}(h)}\lambda _{n}^{\theta _{0}\omega _{2}(h)}\Phi _{n}^{m_{0}(1+\rho
_{h})}(\delta _{\ast })
\end{equation*}%
with%
\begin{equation*}
\xi _{2}(h)=\rho _{h}\xi _{1}(q+2h)=\rho _{h}(q+2h+d+2\theta _{1}+m_{0}(a+b))
\end{equation*}%
and%
\begin{eqnarray*}
\omega _{2}(h) &=&(a+b+\delta _{\ast })m_{0}-\rho _{h}(q+2h+d+2\theta _{1})
\\
&=&\delta _{\ast }m_{0}-\frac{(a+b)m_{0}+q+2d/p_{\ast }}{2h}(q+d+2\theta
_{1})-(q+2d/p_{\ast }).
\end{eqnarray*}%
By our choice of $m_{0}$ we have 
\begin{equation*}
\delta _{\ast }m_{0}>q+2d/p_{\ast }
\end{equation*}%
so, taking $h$ sufficiently large we get $\omega _{2}(h)>0.$ And we also
have $\xi _{2}(h)\leq \xi _{3}:=(a+b)m_{0}+q+\frac{2d}{p_{\ast }}%
+\varepsilon _{\ast }$ and $\rho _{h}\leq \varepsilon _{\ast }.$ So we
finally get%
\begin{equation}
d_{m_{0}}(n)\times \theta (n)^{\rho _{h}}\leq Ct^{-\theta _{0}\xi _{3}}\Phi
_{n}^{m_{0}(1+\varepsilon _{\ast })}(\delta _{\ast }).  \label{R8'}
\end{equation}%
The above inequality guarantees that (\ref{reg10}) holds so that we may use
Lemma \ref{REG}. We take $\eta >\kappa +d$ and, using $(A_{4})$ (see (\ref%
{R7})) we obtain 
\begin{equation*}
\left\vert \mu ^{\eta ,\kappa }\right\vert =\int_{{\mathbb{R}}^{2}}\frac{%
\psi _{\eta }(x)}{\psi _{\kappa }(y)}P_{t}(x,dy)dx\leq C\int_{{\mathbb{R}}}%
\frac{dx}{\psi _{\kappa -\eta }(x)}<\infty .
\end{equation*}%
Then, $A(\delta )<C$ (see (\ref{reg12'})). One also has $B(\varepsilon
)<\infty $ (see (\ref{reg12''})) and finally (see (\ref{reg11})) 
\begin{equation*}
C_{h,n_{\ast }}(\varepsilon )\leq Ct^{-\theta _{0}\xi _{3}}\Phi
_{n}^{m_{0}(1+\varepsilon _{\ast })}(\delta _{\ast }).
\end{equation*}%
We have used here (\ref{R8'}). For large $h$ we also have 
\begin{equation*}
\theta (n)^{\rho _{h}}\leq C(\lambda _{n}t)^{-\theta _{0}((a+b)m_{0}+q+\frac{%
2d}{p_{\ast }})(1+\varepsilon _{\ast })}\Phi _{n}^{\varepsilon _{\ast }}(0).
\end{equation*}

Now (\ref{reg12}) gives (\ref{TR6'}). So \textbf{A }and \textbf{B }are
proved.

\smallskip

\textbf{Step 5: proof of C.} We apply \textbf{B.} with $q$ replaced by $\bar{%
q}=q+1$, so $\Psi _{\eta ,\kappa }p_{t}\in W^{\bar{q},p}({\mathbb{R}}^{d}\times {%
\mathbb{R}}^{d})=W^{\bar{q},p}({\mathbb{R}}^{2d})$. Since $\bar{q}>2d/p$
(here the dimension is $2d$), we can use the Morrey's inequality: for every $%
\alpha $, $\beta $ with $|\alpha |+|\beta |\leq \lfloor \bar q-2d/p\rfloor=q$%
, then $|\partial _{x}^{\alpha }\partial _{y}^{\beta }(\Psi _{\eta
,\kappa }p_{t})(x,y)|\leq C\Vert \Psi _{\eta ,\kappa }p_{t}\Vert _{\bar{q}%
,p} $. By (\ref{TR6'}), one has

\begin{equation*}
\left\vert \partial _{x}^{\alpha }\partial _{y}^{\beta }(\Psi _{\eta ,\kappa
}p_{t})(x,y)\right\vert 
\leq CQ_{\ast} ( \bar q, \bar{\mathfrak{m}} ) 
t^{-\theta_0((a+b)
	\bar{\mathfrak{m}}+\bar{q}+2d/p_{\ast })(1+\varepsilon_{\ast})}
\end{equation*}%
i.e. (using (\ref{NOT3c})), 
\begin{equation*}
\left\vert \partial _{x}^{\alpha }\partial _{y}^{\beta
}p_{t}(x,y)\right\vert
\leq CQ_{\ast} ( \bar q, \bar{\mathfrak{m}} ) 
t^{-\theta_0((a+b)
	\bar{\mathfrak{m}}+\bar{q}+2d/p_{\ast })(1+\varepsilon_{\ast})}\times \frac{1}{\Psi _{\eta ,\kappa
	}(x,y)}.
\end{equation*}%
Now, by a standard calculus, ${\Psi _{\eta ,\kappa }(x,y)}\geq C_{\kappa }%
\frac{\psi _{\kappa }(x-y)}{\psi _{\eta +\kappa }(x)}$ (use that $\psi
_{\kappa }(x-y)\leq C_{\kappa }\psi _{\kappa }(x)\psi _{\kappa
}(-y)=C_{\kappa }\psi _{\kappa }(x)\psi _{\kappa }(y)$), so (\ref{TR6d})
follows. $\square $

\medskip

We are finally ready for the

\smallskip

\textbf{Proof of Theorem \ref{TransferBIS}.} Our assumptions guarantees that 
$P_{t}(x,dy)=p_{t}(x,y)dy$ and $p_{t}$ satisfies (\ref{TR6d}). We take a
cut-off function $F_{R}\in C^{\infty }({\mathbb{R}}^{d})$ such that $%
1_{B_{R}(0)}\leq F_{R}\leq 1_{B_{R+1}(0)}$ ($B_r(0)$ denoting the open ball
centered at $0$ with radius $r$) and we denote $%
p_{t}^{R}(x,y)=F_{R}(x)p_{t}(x,y).$ By (\ref{TR6d}) we know that, for every $%
\kappa \in {\mathbb{N}},\varepsilon >0$ and every $(x,y)\in {\mathbb{R}}%
^{d}\times {\mathbb{R}}^{d}$ one has 
\begin{equation*}
\left\vert \partial _{x}^{\alpha }\partial _{y}^{\beta
}p_{t}^{R}(x,y)\right\vert \leq Ct^{-\theta _{\ast }(\left\vert \alpha
\right\vert +\left\vert \beta \right\vert +\theta _{1})}\times \frac{\psi
_{\eta +\kappa }(x)}{\psi _{\kappa }(x-y)} 
\end{equation*}
%
%
where $\theta _{\ast }=\theta _{0}(1+\frac{a+b}{\delta _{\ast }}%
)(1+\varepsilon )$, $\theta _{1}$ is computed from (\ref{TR6d}) (the precise
value is not important here) and $C$ and $\eta$ both depend on 
$\kappa ,\varepsilon ,\delta_\ast,\left\vert \alpha \right\vert $ and $%
\left\vert \beta \right\vert .$ Since the above left hand side is
identically null when $|x|>R+1$, we can write 
\begin{equation*}
\left\vert \partial _{x}^{\alpha }\partial _{y}^{\beta
}p_{t}^{R}(x,y)\right\vert \leq Ct^{-\theta _{\ast }(\left\vert \alpha
\right\vert +\left\vert \beta \right\vert +\theta _{1})}\times \psi
_{-\kappa }(x,y) 
\end{equation*}
where $C$ is a new constant depending on $R$ as well (we also stress that
here $\psi _{-\kappa }(x,y)=(1+|x|^2+|y|^2)^{-\kappa}$, so the underlying
dimension is $2d$). This allows one to apply Lemma \ref{reg}: for every $%
p\geq 1$, $\kappa\in{\mathbb{N}}$ and $\delta>0$, 
\begin{equation*}
\left\Vert \psi _{\kappa }p_{t}^{R}\right\Vert _{q,p}\leq Ct^{-\theta _{\ast
}(q+\frac{2d}{p_{\ast }}+\delta)},\quad t<2^ {-4/\theta\ast}.
\end{equation*}%
Then by Morrey's Lemma, for every $p>2d$%
\begin{equation*}
\left\Vert \psi _{\kappa }p_{t}^{R}\right\Vert _{q,\infty }\leq \left\Vert
\psi _{\kappa }p_{t}^{R}\right\Vert _{q+1,p}\leq Ct^{-\theta _{\ast }(q+1+%
\frac{2d}{p_{\ast }}+\delta)}\leq Ct^{-\theta _{\ast }(q+2d+\varepsilon )}
\end{equation*}%
the last inequality being true if we take $p$ close to $2d$ and $%
\delta<\varepsilon$. And this gives (\ref{TR6e}). $\square $

%
%


\appendix

\section{Weights}

\label{app:weights}

For $k\in{\mathbb{Z}}$ and $x\in{\mathbb{R}}^d$, we denote%
\begin{equation}
\psi _{k}(x)=(1+\left\vert x\right\vert ^{2})^{k}.  \label{n1}
\end{equation}

\begin{lemma}
\label{Psy1}For every multi-index $\alpha $ there exists a constant $%
C_{\alpha }$ such that 
\begin{equation}
\Big|\partial ^{\alpha }\Big(\frac{1}{\psi _{k}}\Big)\Big\vert \leq \frac{%
C_{\alpha }}{\psi _{k}}.  \label{n2}
\end{equation}%
Moreover, for every $q$ there is a constant $C_{q}\geq 1$ such that for
every $f\in C_{b}^{\infty }({\mathbb{R}}^{d})$%
\begin{equation}
\frac{1}{C_{q}}\sum_{0\leq |\alpha | \leq q}\Big\vert \partial ^{\alpha }%
\Big(\frac{f}{\psi _{k}}\Big)\Big\vert \leq \sum_{0\leq \left\vert \alpha
\right\vert \leq q}\frac{1}{\psi _{k}}\left\vert \partial ^{\alpha
}f\right\vert \leq C_{q}\sum_{0\leq \left\vert \alpha \right\vert \leq q}%
\Big\vert \partial ^{\alpha }\Big(\frac{f}{\psi _{k}}\Big)\Big\vert .
\label{n3}
\end{equation}
\end{lemma}

\textbf{Proof.} One checks by recurrence that 
\begin{equation*}
\partial ^{\alpha }\Big(\frac{1}{\psi _{k}}\Big)=\sum_{q=1}^{\left\vert
\alpha \right\vert }\frac{P_{\alpha ,q}}{\psi _{k+q}}
\end{equation*}%
where $P_{\alpha ,q}$ is a polynomial of order $q.$ And since%
\begin{equation*}
\frac{(1+\left\vert x\right\vert )^{q}}{(1+\left\vert x\right\vert
^{2})^{q+k}}\leq \frac{C}{(1+\left\vert x\right\vert ^{2})^{k}}
\end{equation*}%
the proof (\ref{n2}) is completed. In order to prove (\ref{n3}) we write%
\begin{equation*}
\partial ^{\alpha }\Big(\frac{f}{\psi _{k}}\Big)=\frac{1}{\psi _{k}}\partial
^{\alpha }f+\sum_{\substack{ (\beta ,\gamma )=\alpha  \\ \left\vert \beta
\right\vert \geq 1}}c(\beta ,\gamma )\partial ^{\beta }\Big(\frac{1}{\psi
_{k}}\Big)\partial ^{\gamma }f.
\end{equation*}%
This, together with (\ref{n2}) implies 
\begin{equation*}
\Big\vert \partial ^{\alpha }\Big(\frac{f}{\psi _{k}}\Big)\Big\vert \leq
C\sum_{0\leq \left\vert \gamma \right\vert \leq \left\vert \alpha
\right\vert }\frac{1}{\psi _{k}}\left\vert \partial ^{\gamma }f\right\vert
\end{equation*}%
so the first inequality in (\ref{n3}) is proved. In order to prove the
second inequality we proceed by recurrence on $q$. The inequality is true
for $q=0.$ Suppose that it is true for $q-1.$ Then we write%
\begin{equation*}
\frac{1}{\psi _{k}}\partial ^{\alpha }f=\partial ^{\alpha }\Big(\frac{f}{%
\psi _{k}}\Big)-\sum_{\substack{ (\beta ,\gamma )=\alpha  \\ \left\vert
\beta \right\vert \geq 1}}c(\beta ,\gamma )\partial ^{\beta }\Big(\frac{1}{%
\psi _{k}}\Big)\partial ^{\gamma }f
\end{equation*}%
and we use again (\ref{n2}) in order to obtain%
\begin{equation*}
\frac{1}{\psi _{k}}\left\vert \partial ^{\alpha }f\right\vert \leq \Big\vert %
\partial ^{\alpha }\Big(\frac{f}{\psi _{k}}\Big)\Big\vert +C\sum_{\left\vert
\gamma \right\vert <\left\vert \alpha \right\vert }\frac{1}{\psi _{k}}%
\left\vert \partial ^{\gamma }f\right\vert \leq C\sum_{0\leq \left\vert
\beta \right\vert \leq q}\Big\vert \partial ^{\beta }\Big(\frac{f}{\psi _{k}}%
\Big)\Big\vert
\end{equation*}%
the second inequality being a consequence of the recurrence hypothesis. $%
\square $

\begin{remark}
The assertion is false if we define $\psi _{k}(x)=(1+\left\vert x\right\vert
)^{k}$ because $\partial _{i}\partial _{j}\left\vert x\right\vert =\frac{%
\delta _{i,j}}{\left\vert x\right\vert }-\frac{x_{i}x_{j}}{\left\vert
x\right\vert ^{2}}$ blows up in zero.
\end{remark}

We look now to $\psi _{k}$ itself.

\begin{lemma}
\label{Psy2}For every multi-index $\alpha $ there exists a constant $%
C_{\alpha }$ such that 
\begin{equation}
\left\vert \partial ^{\alpha }\psi _{k}\right\vert \leq C_{\alpha }\psi _{k}.
\label{n4}
\end{equation}%
Moreover, for every $q$ there is a constant $C_{q}\geq 1$ such that for
every $f\in C_{b}^{\infty }({\mathbb{R}}^{d})$%
\begin{equation}
\frac{1}{C_{q}}\sum_{0\leq \left\vert \alpha \right\vert \leq q}\left\vert
\partial ^{\alpha }(\psi _{k}f)\right\vert \leq \sum_{0\leq \left\vert
\alpha \right\vert \leq q}\psi _{k}\left\vert \partial ^{\alpha
}f\right\vert \leq C_{q}\sum_{0\leq \left\vert \alpha \right\vert \leq
q}\left\vert \partial ^{\alpha }(\psi _{k}f)\right\vert .  \label{n5}
\end{equation}
\end{lemma}

\textbf{Proof.} One proves by recurrence that, if $\left\vert \alpha
\right\vert \geq 1$ then $\partial ^{\alpha }\psi
_{k}=\sum_{q=1}^{\left\vert \alpha \right\vert }\psi _{k-q}P_{q}$ with $%
P_{q} $ a polynomial of order $q.$ Since $1+\left\vert x\right\vert \leq
2(1+\left\vert x\right\vert ^{2})$ it follows that $\left\vert
P_{q}\right\vert \leq C\psi _{q}$ and (\ref{n4}) follows. Now we write 
\begin{equation*}
\psi _{k}\partial ^{\alpha }f=\partial ^{\alpha }(\psi _{k}f)-\sum 
_{\substack{ (\beta ,\gamma )=\alpha  \\ \left\vert \beta \right\vert \geq 1 
}}c(\beta ,\gamma )\partial ^{\beta }\psi _{k}\partial ^{\gamma }f
\end{equation*}%
and the same arguments as in the proof of (\ref{n3}) give (\ref{n5}).

\section{Semigroup estimates}

\label{app:semi}

We consider a semigroup $(P_t)_{t\geq 0}$ on $C^{\infty }({\mathbb{R}}^{d})$
such that $P_{t}f(x)=\int f(y)P_{t}(x,dy)$ where $P_{t}(x,dy)$ is a
probability transition kernel and we denote by $P_{t}^{\ast }$ its formal
adjoint.

\begin{assumption}
\label{B1B2} 
There exists $Q\geq 1$ such that for every $t\leq T$ and every $f\in
C^{\infty }({\mathbb{R}}^{d})$%
\begin{equation}
\left\Vert P_{t}f\right\Vert _{1}\leq Q\left\Vert f\right\Vert _{1}.
\label{A31}
\end{equation}
Moreover, 
for every $k\in {\mathbb{N}}$ there exists $K_{k}\geq 1$ such that for every 
$x\in {\mathbb{R}}^{d}$ and $t\leq T$
\begin{equation}
\left\vert P_{t}(\psi _{k})(x)\right\vert \leq K_{k}\psi _{k}(x).
\label{A32}
\end{equation}
\end{assumption}

\begin{lemma}
Under Assumption \ref{B1B2}, for every $t\leq T$ one has 
\begin{equation}
\left\Vert \psi _{k}P_{t}^{\ast }(f/\psi _{k})\right\Vert _{p}\leq
K_{kp}^{1/p}Q^{1/p_{\ast }}\left\Vert f\right\Vert _{p}.  \label{A34}
\end{equation}
\end{lemma}

\textbf{Proof.} Using H\"{o}lder's inequality, the identity $\psi
_{k}^{p}=\psi _{kp},$ and (\ref{A32})%
\begin{equation*}
\left\vert P_{t}(\psi _{k}g)(x)\right\vert \leq \left\vert P_{t}(\psi
_{k}^{p})(x)\right\vert ^{1/p}\left\vert P_{t}(\left\vert g\right\vert
^{p_{\ast }})(x)\right\vert ^{1/p_{\ast }}\leq K_{kp}^{1/p}\psi
_{k}(x)\left\vert P_{t}(\left\vert g\right\vert ^{p_{\ast }})(x)\right\vert
^{1/p_{\ast }}.
\end{equation*}%
Then, using (\ref{A31}) 
\begin{eqnarray*}
\left\Vert \frac{1}{\psi _{k}}P_{t}(\psi _{k}g)\right\Vert _{p_{\ast }}
&\leq &K_{kp}^{1/p}\left\Vert \left\vert P_{t}(\left\vert g\right\vert
^{p_{\ast }})\right\vert ^{1/p_{\ast }}\right\Vert _{p_{\ast
}}=K_{kp}^{1/p}(\left\Vert P_{t}(\left\vert g\right\vert ^{p_{\ast
}})\right\Vert _{1})^{1/p_{\ast }} \\
&\leq &K_{kp}^{1/p}Q^{1/p_{\ast }}(\left\Vert \left\vert g\right\vert
^{p_{\ast }}\right\Vert _{1})^{1/p_{\ast }}=K_{kp}^{1/p}Q^{1/p_{\ast
}}\left\Vert g\right\Vert _{p_{\ast }}.
\end{eqnarray*}%
Using H\"{o}lder's inequality first and the above inequality we obtain%
\begin{eqnarray*}
\left\vert \left\langle g,\psi _{k}P_{t}^{\ast }(f/\psi _{k})\right\rangle
\right\vert &=&\left\vert \left\langle \frac{1}{\psi _{k}}P_{t}(g\psi
_{k}),f\right\rangle \right\vert \leq \left\Vert f\right\Vert _{p}\left\Vert 
\frac{1}{\psi _{k}}P_{t}(g\psi _{k})\right\Vert _{p_{\ast }} \\
&\leq &K_{kp}^{1/p}Q^{1/p_{\ast }}\left\Vert g\right\Vert _{p_{\ast
}}\left\Vert f\right\Vert _{p}.
\end{eqnarray*}%
$\square $

\smallskip

We consider also the following hypothesis.

\begin{assumption}
\label{B3} There exists $\rho >1$ such that 
for every $q\in {\mathbb{N}}$ there exists $D_{(q)}^{\ast }(\rho )\geq 1$
such that for every $x\in {\mathbb{R}}^{d}$ and $t\leq T$
\begin{equation}
\sum_{\left\vert \alpha \right\vert \leq q}\left\vert \partial ^{\alpha
}P_{t}^{\ast }f(x)\right\vert \leq D_{(q)}^{\ast }(\rho )\sum_{\left\vert
\alpha \right\vert \leq q}(P_{t}^{\ast }(\left\vert \partial ^{\alpha
}f\right\vert ^{\rho })(x))^{1/\rho }.  \label{A36}
\end{equation}
\end{assumption}

\begin{proposition}
\label{A2}Suppose that Assumption \ref{B1B2} and \ref{B3} 
hold. Then for every $k,q\in {\mathbb{N}}$ and $p>\rho $ there exists a
universal constant $C$ (depending on $k$ and $q$ only) such that for every $t\leq T$
\begin{equation}
\left\Vert \psi _{k}P_{t}^{\ast }(f/\psi _{k})\right\Vert _{q,p}\leq
CK_{k p}^{1/p}Q^{(p-\rho )/\rho p}D_{(q)}^{\ast }(\rho )\left\Vert
f\right\Vert _{q,p}.  \label{A38}
\end{equation}
\end{proposition}

\textbf{Proof.} We will prove (\ref{A38}). Let $\alpha $ with $\left\vert
\alpha \right\vert \leq q.$ By (\ref{A36}) 
\begin{eqnarray*}
\left\vert \partial ^{\alpha }(\psi _{k}P_{t}^{\ast }(f/\psi
_{k})(x))\right\vert &\leq &C\psi _{k}(x)\sum_{\left\vert \gamma \right\vert
\leq q}\left\vert \partial ^{\gamma }(P_{t}^{\ast }(f/\psi
_{k})(x))\right\vert \\
&\leq &CD_{(q)}^{\ast }(\rho )\psi _{k}(x)\sum_{\left\vert \beta \right\vert
\leq q}(P_{t}^{\ast }(\left\vert \partial ^{\beta }(f/\psi _{k})\right\vert
^{\rho })(x))^{1/\rho } \\
&=&CD_{(q)}^{\ast }(\rho )\sum_{\left\vert \beta \right\vert \leq q}(\psi
_{\rho k}(x)P_{t}^{\ast }(\left\vert \partial ^{\beta }(f/\psi
_{k})\right\vert ^{\rho })(x))^{1/\rho } \\
&=&CD_{(q)}^{\ast }(\rho )\sum_{\left\vert \beta \right\vert \leq q}(\psi
_{\rho k}(x)P_{t}^{\ast }(g_\beta/\psi _{\rho k})(x))^{1/\rho }
\end{eqnarray*}%
with%
\begin{equation*}
g_\beta(x)=\psi _{\rho k}(x)\left\vert \partial ^{\beta }(f/\psi
_{k})(x)\right\vert ^{\rho }=\left\vert \psi _{k}(x)\partial ^{\beta
}(f/\psi _{k})(x)\right\vert ^{\rho }.
\end{equation*}%
Taking $p>\rho$ and using (\ref{A34}) 
\begin{equation*}
\left\Vert (\psi _{\rho k}P_{t}^{\ast }(g_\beta/\psi _{\rho k}))^{1/\rho
}\right\Vert _{p}=\left\Vert \psi _{\rho k}P_{t}^{\ast }(g_\beta/\psi _{\rho
k})\right\Vert _{p/\rho }^{1/\rho }\leq K_{k p}^{1/p}Q^{(p-\rho )/\rho
p}\left\Vert g_\beta\right\Vert _{p/\rho }^{1/\rho }.
\end{equation*}%
And we have%
\begin{equation*}
\left\Vert g_\beta\right\Vert _{p/\rho }^{1/\rho }=(\int \left\vert \psi
_{k}(x)\partial ^{\beta }(f/\psi _{k})(x)\right\vert ^{p}dx)^{1/p}\leq
C\sum_{\left\vert \gamma \right\vert \leq q}(\int \left\vert \partial
^{\gamma }f(x)\right\vert ^{p}dx)^{1/p}=C\left\Vert f\right\Vert _{q,p}.
\end{equation*}%
We conclude that 
\begin{equation*}
\left\Vert \psi _{k}P_{t}^{\ast }(f/\psi _{k})\right\Vert _{q,p}\leq
CK_{k p}^{1/p}Q^{(p-\rho )/\rho p}D_{(q)}^{\ast }(\rho )\left\Vert
f\right\Vert _{q,p}.
\end{equation*}%
$\square $

\end{document}